\documentclass[12pt]{article}
\usepackage{amsmath,amsfonts,amssymb}
\usepackage[all]{xy}
\advance \textheight by \topskip
\makeatletter
 \def\theequation{\thesection.\arabic{equation}}
\makeatother    

\newtheorem{theorem}{Theorem}[section]

\newtheorem{corollary}[theorem]{Corollary}

\newtheorem{definition}[theorem]{Definition}
\newtheorem{example}[theorem]{Example}

\newtheorem{lemma}[theorem]{Lemma}

\newtheorem{proposition}[theorem]{Proposition}
\newtheorem{remark}[theorem]{Remark}

\def\NN{\mathbb N}
\def\ZZ{\mathbb Z}
\def\CC{\mathbb C}
\def\RR{\mathbb R}
\newcommand{\beq}{\begin{equation}}
\newcommand{\eeq}{\end{equation}}
\newcommand{\beqa}{\begin{eqnarray}}
\newcommand{\eeqa}{\end{eqnarray}}
\newcommand{\g}{{\mathfrak g}}
\newcommand{\n}{{\mathfrak n}}
\newcommand{\rr}{{\mathfrak r}}

\newcommand{\s}{{\mathfrak{sl}(2)}}
\newcommand{\A}{{\mathrm {Aut}}}

\newcommand{\noi}{\noindent}

\def\>{\rangle}
\def\<{\langle}

\begin{document}

\title{
{\bf Finite-dimensional Lie subalgebras of the Weyl algebra}}

\author{
{\sf M. Rausch de Traubenberg}\thanks{e-mail:
rausch@lpt1.u-strasbg.fr}$\,\,$${}^{a},$
{\sf M. J. Slupinski}\thanks{e-mail:
slupins@math.u-strasbg.fr}$\,\,$${}^{b}$
and
{\sf A. Tanas\u{a}}\thanks{e-mail:
atanasa@lpt1.u-strasbg.fr}$\,\,$$^{c,a}$\\
\\
{\small ${}^{a}${\it
Laboratoire de Physique Th\'eorique, CNRS UMR  7085,
Universit\'e Louis Pasteur}}\\
{\small {\it  3 rue de
l'Universit\'e, 67084 Strasbourg Cedex, France}}\\
{\small ${}^{b}${\it
Institut de Recherches en Math\'emathique Avanc\'ee,
Universit\'e Louis Pasteur and CNRS}}\\
{\small {\it 7 rue R. Descartes, 67084 Strasbourg Cedex, France.}}  \\ 
{\small ${}^{c}${\it Laboratoire de Math\'ematiques et Applications,
Universit\'e de Haute Alsace,}} \\
{\small {\it Facult\'e des Sciences et Techniques,
4 rue des Fr\`eres Lumi\`ere, F-68093 Mulhouse, France.}}}   

\date{}

\maketitle
\vskip-1.5cm

\vspace{2truecm}

\begin{abstract}
\noindent
We classify up to isomorphism all finite-dimensional Lie algebras that can be realised as Lie subalgebras of the complex Weyl algebra $A_1$. The list we obtain turns out to be discrete and for example, the only non-solvable Lie algebras with this property are:
 $\s$, $\s\times\CC$ and  $\s\ltimes{\cal H}_3$. 
We then give  several different characterisations, normal forms and  isotropy groups for the action of  $\A (A_1)\times \A (\s)$ on a particular class of realisations of $\s$ in $A_1$.
\end{abstract}

\vspace{2cm}

\newpage
\section{Introduction}
\renewcommand{\theequation}{1.\arabic{equation}}   
\setcounter{equation}{0}

The  Weyl algebra $A_1$ is the complex associative algebra generated by elements $p$ and $q$ satisfying the relation $pq-qp=1$. It is well known that the Lie algebras $\s$, $\s\times\CC$ and $\s\ltimes{\cal H}_3$ (where ${\cal H}_3$ denotes the three-dimensional Heisenberg algebra) 
can be realised as Lie subalgebras of $A_1$. In \cite{italienii},
A. Simoni and F. Zaccaria proved that the only complex semi-simple Lie algebra that can be realised in $A_1$ is $\s$ and a remarkable property of  realisations of $\s$ in $A_1$  was proved by A. Joseph in \cite{joseph}, where he showed that the spectrum of the realisation in $A_1$ of  suitably normalised semi-simple elements of $\s$ is either $\ZZ$ or $2\ZZ$.
In \cite{japonezu}, J. Igusa gave a necessary condition for two elements of $A_1$ to generate an infinite-dimensional Lie subalgebra.

In this article we find all complex finite-dimensional Lie algebras which can be realised in $A_1$. If $\g$ is a complex finite-dimensional Lie algebra and
 $A_1^{\mathfrak{g}}$ 
denotes the set of injective Lie algebra homomorphisms  from $\g$ to $A_1$, we prove (Theorems \ref{mareata}, \ref{ln}, \ref{solvable22} and \ref{solvable2})

\begin{theorem}
\label{intro-1}
Let $\g$ be a complex finite-dimensional non-abelian Lie algebra. Then $A_1^\g\ne\emptyset$ iff $\g$ is isomorphic to one of the following:
\beqa
\begin{array}{llll}
1) & \s, & 4) & {\cal L}_n\ (n\ge 2),\\
2) & \s \times \CC, & 5) & \tilde {\cal L}_n \ (n\ge 2) , \\
3) & \s\ltimes\mathcal H_3, \ \ \ \ \ \ \ \ \  & 6) & \rr (i_1,\dots ,i_n) \ \ (\mbox{$i_1<\ldots <i_n$ are positive integers}). 
\end{array}\nonumber
\eeqa
\end{theorem}
\noi
Here, ${\cal L}_n$ is a  nilpotent, in fact filiform, Lie algebra, $ \tilde {\cal L}_n$ is isomorphic to a  semi-direct product $\CC \ltimes {\cal L}_n$   and 
$\rr (i_1, \dots i_n)$
is isomorphic to a  semi-direct product  $\CC \ltimes \CC^n$. For the precise definitions of these Lie algebras see section $3$.
Note that only a finite number  of non-solvable Lie algebras  and only a discrete family of solvable Lie algebras appear in the list of Theorem \ref{intro-1}.
Since all derivations of $A_1$ are inner, this theorem also leads to the classification of all finite-dimensional Lie algebras which can be realised in Der$(A_1)$ (see Theorem \ref{ptintro}) and we give explicit examples of subgroups of $\A (A_1)$ which exponentiate them (see section $4.4$).

In the second part of the paper we study a particular family of realisations of $\s$ in $A_1$. 
The group $\A (A_1)\times\A(\s)$
acts naturally on $A_1^\s$ and we give several characterisations of  the orbit of 
 ${\cal N} = \{ f_I,\ f_{II}^{b} \ : \ b\in{\CC}\}\subseteq A_1^\s$ 
where
\beqa
\begin{array}{ll}
 f_I(e_+)=-\frac 12 q^2 &
f_{II}^{b}(e_+)= (b+pq)q \\
f_I(e_-)=\frac 12 p^2 &
f_{II}^{b}(e_-)=-p \\
f_I(e_0)=\frac 12 (pq+qp) &
f_{II}^{b}(e_0)=2pq+b
\end{array}
\eeqa
\noi
( $e_+,e_-,e_0$ is the standard basis of $\s$). The realisations $f_{II}^{b}$ were first introduced in this context by A. Joseph in \cite{joseph} and $f_I$
is the natural  embedding of $\s$ in $A_1$. We  prove that no two elements of  $\cal N$ are in the same orbit and the
second main result of this paper 
is (Theorem \ref{TEOREMA})

\begin{theorem}
\label{intro-2}
Let $f\in A_1^{\s}$. 
Then the following statements are equivalent:
\begin{enumerate}
\item[] (i) There exists $\gamma\in \A (A_1)\times \A (\s)$ such that  $\gamma.f\in \cal N$;
\item[] (ii) There exists $z\in\s\setminus\{ 0 \}$  such that ${\rm ad}(f(z))$ has a strictly semi-simple eigenvector;
\item[] (iii)  There  exists $z\in\s\setminus\{ 0 \}$  such that $f(z)$ is strictly semi-simple;
\item[] (iv) There exists $z\in\s\setminus\{ 0 \}$  such that ${\rm ad}(f(z))$ has a strictly nilpotent eigenvector;
\item[] (v) There exists  $z\in\s\setminus\{ 0 \}$  such that  $f(z)$ is strictly nilpotent;
\item[] (vi) There exists $z\in\s\setminus\{ 0 \}$ such that $\mathrm{ad}(f(z))$ can be exponentiated in $\A (A_1)$.
\end{enumerate}
\end{theorem}
\noi
The terms strictly nilpotent and strictly semi-simple for non-zero elements in  $A_1$ were defined by J. Dixmier in \cite{Dixmier} (see also subsection $2.2$ of this paper) and for the precise definition of exponentiation in this context, see subsection $2.3$.
We further show that the isotropy of $f_I$ is isomorphic to $SL(2,\CC)$ and that the isotropy of $f_{II}^{b}$ is isomorphic to a Borel subgroup of $SO(3,\CC)$.
Finally, for the sake of completeness, we give explicit formulae for a realisation of $\s$ in $A_1$ which does not satisfy any of the criteria of  Theorem \ref{intro-2}. (see also \cite{joseph}).

\medskip

The paper is organised as follows. In section $2$ we recall the basic properties of the Weyl algebra and in particular the Dixmier partition which is essential to this article. In section $3$ we give some examples of Lie algebras which can be realised as Lie subalgebras of $A_1$ and in section $4$ we obtain the classification of all finite-dimensional Lie algebras with this property. Sections $5$ and $6$ are devoted to the study of an explicit  family  of realisations of $\s$ and its orbit under the action of the group $\A (A_1)\times \A (\s)$.

\section{Properties of the Weyl algebra and the Dixmier partition}
\renewcommand{\theequation}{2.\arabic{equation}}   
\setcounter{equation}{0}

In this section we  give the  basic properties of 
  the Weyl algebra $A_1$ and its  Dixmier partition. In particular we give a simple characterisation of  the set of elements $Z\in A_1$ such that $\mathrm{ad} (Z)$ can be exponentiated.

\subsection{Basic properties}

\begin{definition}
\label{Weyl}
The \textit{ Weyl algebra} $A_1$ is the complex associative algebra generated by elements $p$ and $q$ subject to the relation $pq-qp=1$.
\end{definition}

There is a natural action of $A_1$ on ${\CC}[x]$ defined by
\beq
\label{pq}
p\cdot P(x) = P'(x), \ \ q\cdot P(x) =x P(x),\ \forall P\in {\CC}[x]
\eeq
and this establishes an isomorphism of $A_1$ with
the algebra of polynomial coefficient differential operators acting on 
 ${\CC}[x]$. 
We will refer to this representation as the \textit{standard representation of $A_1$}.

\medskip

\noindent
\textbf{Properties of $A_1$} (see \cite{Dixmier}, for example):
\begin{enumerate}
\item [{\bf  P1}] The elements $\{p^iq^j : i,j\in \mathbb N\}$ constitute a basis
of $A_1$ and the centre of $A_1$  is $\CC$.
 \item [{\bf P2}] The linear subspace $W_1= <p,q>$ spanned by $p$ and $q$ has a
unique symplectic structure $\omega$ such that $\omega(p,q)=1$ and 
 the group 
$SL(W_1)$ of symplectic transformations of ($W_1,\omega$) is isomorphic to
$SL(2,\mathbb C)$.
It is well
 known that the inclusion of $W_1$ in $A_1$ extends to an algebra isomorphism from the quotient of the tensor
 algebra $T(W_1)$ by the two-sided ideal generated by $v_1\otimes v_2
 -v_2\otimes v_1 - \omega (v_1,v_2)1$ to $A_1$ and, since this ideal
 is stable under its action, $SL(W_1)$ acts naturally on $A_1$.
The map 
$\delta : S(W_1) \to A_1$ given on $S^n(W_1)$ by 
$$\delta(v_1\odot\ldots\odot v_n)=\sum_{\sigma\in S_n} \frac {1}{n!}
 v_{\sigma(1)}\dots v_{\sigma(n)}$$
is an $SL(W_1)$-equivariant  {\it vector space}
isomorphism. 
If we set $W_n=\delta (S_n (W_1))$  one can show that
 $A_1=\bigoplus_{n\in\NN} W_n$ and that $[W_i,W_j]\subseteq
 W_{i+j-2}$. In particular,
\beqa
\label{W2}
[W_2, W_2]\subseteq W_2,\ [W_2,W_i]\subseteq W_i
\eeqa
so that $W_2$ is a Lie algebra and $W_i$ is a representation of
$W_2$. 
This action of $W_2$ on $W_1$ establishes an $SL(W_1)$-equivariant Lie algebra isomorphism
$W_2\cong {\mathfrak {sl}}(W_1)$. 
\item [{\bf  P3}] The algebra $A_1$ satifies the commutative centraliser condition
 (\textit{ccc}): the centraliser $C(x)$ of any element  $x\in A_1\setminus \CC$ is a commutative subalgebra (see \cite{Amitsur} and \cite{Dixmier}). 
\item [{\bf  P3}${}^\prime$]  If $x,y\in A_1\setminus \CC$ then
$$ C(x)\cap C(y)=\left\{ 
\begin{array}{ll}
\CC & \mbox{ if } xy-yx\ne 0 \\
C(x) & \mbox{ if } xy-yx=0,
\end{array} \right.
$$
(see Corollary $4.7$ of \cite{Dixmier}).
\item[{\bf  P4}]  Two elements $p',q'\in A_1$ satisfying $[p',q']=1$ uniquely define an algebra
homomorphism from $A_1$ to itself and conversely, given an homomorphism $\alpha: A_1\to A_1$ we have $[\alpha (p),\alpha (q)]=1$. J. Dixmier
conjectured in 1968 that an algebra homomorphism of $A_1$ is  invertible
and thus in fact an automorphism (see \cite{Dixmier}). This
conjecture is still undecided.
\item [{\bf  P5}] For $n\in \mathbb N$, the map
 $\mathrm{ad}(p^{n+1}):A_1\to A_1$ given by
$\mathrm{ad}({p^{n+1}})(a)=[p^{n+1},a]$ 
is locally  nilpotent, \textit{i.e.}, for each $a \in A_1$, there exists an $N\in \mathbb N$
(depending on $a$)
such that $\mathrm{ad}^N({p^{n+1}})(a)=0$. For $\lambda\in\CC$ one can then define $\Phi_{n,\lambda}:A_1\to
A_1$ by
$\Phi_{n,\lambda}(a)=\sum_{k=0}^N
\frac{(\frac{\lambda}{n+1}ad({p^{n+1}}))^k}{k!}(a)=
\exp(\frac{\lambda}{n+1}\mathrm{ad}({p^{n+1}}))(a)$ and this is the unique automorphism of $A_1$ such that
\beqa
\label{auto}
\Phi_{n,\lambda}(p)&=&p, \ \ \Phi_{n,\lambda}(q)=q+\lambda p^n. 
\eeqa
One defines $\Phi'_{n,\lambda}=\mathrm{exp}(-\frac{\lambda}{n+1}\mathrm{ad}({q^{n+1}}))$ similarly and shows that it is the unique automorphism of $A_1$ such that
\beqa
\label{auto'}
\Phi'_{n,\lambda}(q)&=&q, \ \ \Phi'_{n,\lambda}(p)=p+\lambda q^n.
\eeqa
The group of automorphisms of $A_1$ is generated by the
$\Phi_{n,\lambda}$ and the
$\Phi'_{n,\lambda}$ (see \cite{Dixmier}).
\end{enumerate}

\subsection{The Dixmier partition}

 Let $x\in A_1$. Set
\beqa
N(x)&=&\big\{ y\in A_1 :\ \mathrm{ad}^m(x) (y)=0, \mbox{ for some positive
integer } m \big\} \nonumber \\
C(x)&=&\big\{ y\in A_1:\  \mathrm{ad}(x)(y)=0 \big\} \nonumber \\
D(x)&=&\Big<y\in A_1 : \mathrm{ad}(x) (y)=\lambda y  \mbox{ for some }
\lambda\in\CC \Big>.\nonumber 
\eeqa
It is immediate that $N(x)\cap D(x)=C(x)$ and Dixmier showed that  for
all $x\in A_1$, either $C(x)=N(x)$ or $C(x)=D(x)$.
As a consequence he proved (see \cite{Dixmier})

\begin{theorem}
\label{th-Dixmier}
(Dixmier partition) The set $A_1\setminus\CC$ is a disjoint union of the following
non-empty subsets.
\beqa
\Delta_1 &=& \big\{ x\in A_1\setminus\CC:\ D(x)=C(x), \ N(x)\ne C(x),\ N(x)=A_1  \big\}
\nonumber\\
\Delta_2 &=& \big\{ x\in A_1\setminus\CC:\ D(x)=C(x), \ N(x)\ne C(x),\ N(x)\ne A_1 \big\} \nonumber\\
\Delta_3 &=& \big\{ x\in A_1\setminus\CC:\ D(x)\ne C(x),\ N(x)=C(x), \ D(x)=A_1\big\}
\nonumber\\
\Delta_4 &=& \big\{ x\in A_1\setminus\CC:\ D(x)\ne C(x),\ N(x)=C(x), \ D(x)\ne A_1\big\} \nonumber\\
\Delta_5 &=& \big\{ x\in A_1\setminus\CC:\ D(x)=C(x),\ N(x)=C(x),\  C(x)\ne A_1\big\} \nonumber
\eeqa
\end{theorem}
\noi
Note that this partition is stable under the action of $\A (A_1)$ and
multiplication by a non-zero scalar.

Elements of $\Delta_1 \cup \Delta_2$ (resp. of $\Delta_1$) are said to
be nilpotent (resp. strictly nilpotent) and elements of $\Delta_3 \cup \Delta_4$ (resp. of $\Delta_3$) are said to
be semi-simple (resp. strictly semi-simple).
In fact, $x\in \Delta_1$ iff there exists an automorphism $\alpha$ of $A_1$ such
that $\alpha (x)$ is a polynomial in $p$ (Theorem $9.1$ of
\cite{Dixmier}) and   $x\in \Delta_3$ iff there exists an automorphism $\alpha$ of $A_1$ such
that $\alpha (x)=\mu pq + \nu$, for some $\mu\in\CC^*$ and $\nu\in\CC$ (Theorem $9.2$ of
\cite{Dixmier}).  

Recall (see {\bf P2}) that $W_2$ is a Lie algebra isomorphic to $\s$ and therefore elements
 of $W_2$ are either semi-simple or nilpotent in the Lie algebra sense.

\begin{proposition}
\label{iar-dixmier}
(\cite{Dixmier}, Lemme $8.6$). Let $x\in A_1\setminus \CC$ be of the form
$ \alpha + w_1 + w_2 $
where $\alpha\in \CC$, $ w_1\in W_1$ and $ w_2\in W_2$.\\
(i) If $w_2$ is nilpotent in the Lie algebra sense, then $x\in \Delta_1$.\\
(ii) If $w_2$ is semi-simple in the Lie algebra sense, then $x\in \Delta_3$.
\end{proposition}

\subsection{Characterisation of $\Delta_1\cup \Delta_3$  in terms of exponentiation}

\begin{definition}
\label{definitie}
Let $Z\in A_1$. One says that ${\mathrm{ad}}(Z)$ can be exponentiated if there exists a
group homomorphism $\Phi:\CC\to \A(A_1)$ such that
\begin{enumerate}
\item for all $a\in A_1$, the vector space  
$V_a=<\Phi(t)(a):\ t\in \CC>$ is finite-dimensional;
\item  $\Phi_a: \CC\to V_a$ is holomorphic and  $\left.\frac{d}{dt}\right\vert_0 \Phi_a (t) = [Z,a]$ where $\Phi_a(t)=\Phi (t)(a)$.
(Since $A_1$ is infinite-dimensional, we impose $1$ so that $2$ makes sense).
\end{enumerate}
\end{definition}

\begin{example}
\label{ex-integrabil} 
If $Z\in A_1$ is such that
$\mathrm{ad}(Z)$ is locally nilpotent 
then $\mathrm{ad}(Z)$ can be exponentiated in this sense (cf. {\bf P5} when $Z$ is a polynomial in $p$).
If $Z=pq+\alpha$, $\mathrm{ad}(Z)$ can be  exponentiated by the group homomorphism 
 $\Psi:\CC\to\A (A_1)$ given on the canonical basis $<p^iq^j:\ i,j\in\NN> $ by $\Psi (t)(p^iq^j)=e^{t(j-i)}p^iq^j$.
\end{example}

\begin{lemma}
\label{exponentiere}
Suppose ${\mathrm{ad}}(Z)$ can be 
exponentiated in the above sense. 
\begin{enumerate}
\item For all $a\in A_1$, the finite-dimensional vector space $V_a$ is stable under the action of $\mathrm{ad}(Z)$.
\item $\Phi (t)(a)=e^{t\mathrm{ad}(Z)\vert_{V_a}}(a)$ for all $a\in A_1$
($e^{t \mathrm{ad}(Z)\vert_{V_a}}$ is well defined by $1$).
\end{enumerate}
\end{lemma}
{\bf Proof:} Fix $a\in A_1$ and $t_0\in \CC$. Then 
$$\Phi (t+t_0)(a) = \Phi (t) \Phi (t_0) (a)$$
and both sides
are in the finite-dimensional vector space $V_a$ so we can differentiate with respect to $t$. This gives
\beqa
\label{muie69}
\left.\frac{d}{dt}\right\vert_{0} \Phi (t+t_0) (a)=\left.\frac{d}{dt}\right\vert_{0} \Phi (t)
\Phi (t_0) (a)=[Z, \Phi (t_0) (a)].
\eeqa
The LHS is in $V_a$ and hence 
$[Z, \Phi (t_0) (a)]$ also. This proves part $1$.

Part $2$ follows from the fact that the curves $t\mapsto \Phi (t)(a)$ and $t\mapsto e^{t\mathrm{ad}(Z)\vert_{V_a}}(a)$
are contained in the finite-dimensional vector space $V_a$ and are solutions of the same first order differential 
equation
$$ \frac{d}{dt}  \gamma(t)  =[Z, \gamma (t) ]$$
with the same initial condition $\gamma(0)=a.$ QED

\begin{proposition}
\label{caracterizare}
Let $Z\in A_1\setminus\{0\}.$ Then 
${\mathrm{ad}}(Z)$ can be  exponentiated iff $Z\in\Delta_1\cup\Delta_3$. 
\end{proposition}
{\bf Proof}: ($\Rightarrow$) : Suppose  that
$Z\notin\Delta_1\cup\Delta_3$. 
Set $F(Z)=\Big\{
a\in A_1:\ \mathrm{dim} (<\mathrm{ad}^n (Z) (a), n\in\NN>)<\infty\Big\}$.
Then by Corollary $6.6$ of \cite{Dixmier},
$F(Z)=D(Z)\cup N(Z)$
and  so by Theorem
\ref{th-Dixmier}, $F(Z)\ne A_1$.
Let $a\in
A_1\setminus F(Z)$.
By hypothesis, $V_a$ is finite-dimensional, stable under ad$(Z)$ and
contains $a$. Hence $<\mathrm{ad}^n (Z) (a):\ n\in\NN>\;\subseteq V_a$ is
also finite-dimensional, which is a contradiction. 

\noi
($\Leftarrow$) : If  $Z\in\Delta_1\cup\Delta_3$ then up to an
automorphism of $A_1$, $Z$ is equal to $pq+\alpha$ 
or to a polynomial in $p$. The result follows from Example
\ref{ex-integrabil}. QED

\section{Examples of Lie subalgebras of the Weyl algebra}
\renewcommand{\theequation}{3.\arabic{equation}}   
\setcounter{equation}{0}

In this section we give  examples of Lie algebras
 which can be realised as Lie subalgebras of $A_1$ (see also \cite{Dixmier}, \cite{japonezu}, \cite{joseph} and \cite{italienii}).

\begin{definition}
Let $\mathfrak{g}$ be a  complex Lie algebra.
$$ A_1^{\mathfrak{g}}=\{f\in{\rm Hom}(\mathfrak{g},A_1): \ f
\mbox{ is injective and } f([a,b])=f(a)f(b)-f(b)f(a) \quad \forall a,b\in
A_1\}.$$
If $A_1^\g\ne \emptyset$ we will say that the Lie algebra $\g$ can be realised as a Lie subalgebra of $A_1$ and 
an element of $A_1^\g$ will be called a realisation of $\g$ in $A_1$.
\end{definition}

\begin{remark}
\label{simplificari}
Let $f:\g\to A_1$ be a realisation of $\g$ and $x\in \g\setminus Z_\g$.
Then $f(x)$ is semi-simple (resp. nilpotent) in the $A_1$ sense if $x$ is semi-simple (resp. nilpotent) in the Lie algebra sense.
For example, when $\mathrm{ad}(x)$ is diagonalisable, there exist $y\in\g$ and $\lambda\in\CC^*$ such that  $[x,y]=\lambda y$;  hence  $f(y)\in D(f(x))$, $f(y)\notin C(f(x))$ and $f(x)$ is semi-simple by Theorem \ref{th-Dixmier}.
\end{remark}

Let us now give some examples of Lie algebras $\g$ for which $A_1^\g\ne\emptyset$.

$  $
\begin{enumerate}
\label{initiatic}
\item [{\bf E1}] We saw in {\bf P2} that $W_2$ is a Lie subalgebra of $A_1$ isomorphic to $\s$. The standard basis 
\begin{eqnarray}
\label{sl2-initiatic}
X = -\frac 12 q^2, \
Y = \frac 12  p^2, \
H = \frac 12 (pq+qp)= pq-\frac 12,  
\end{eqnarray}
satisfies the commutation relations $[X,Y]=H$, $[H,X]=2X$ and $[H,Y]=-2Y$.

\item [{\bf E2}] The  elements $1, X, Y, H$ 
span $W_0\oplus W_2$ which is a Lie subalgebra of $A_1$ isomorphic to the direct
product $\s\times\CC$. 
\item [{\bf E3}] The  elements $1,p,q$  span $W_0\oplus W_1$ which is a Lie subalgebra
isomorphic to the three dimensional Heisenberg algebra 
${\cal H}_3$. 
\item [{\bf E4}] The elements  $1,p,q,X,Y,H$  span $W_0\oplus W_1\oplus W_2$ which is a Lie
subalgebra of $A_1$ isomorphic to a semi-direct product
$\s\ltimes{\cal H}_3$.  
\item [{\bf E5}]
The associative subalgebra of $A_1$ generated by $p$ is infinite-dimensional abelian and therefore  any finite-dimensional abelian Lie algebra can be realised in $A_1$.
 \item [{\bf E6}]
The $(n+1)$ elements $q, 1, p, \ldots, p^{n-1}$
 span a non-abelian nilpotent Lie subalgebra isomorphic to
the filiform Lie algebra ${\cal L}_n$ (see \cite{muierea}).  If we set $X_0=-q,
X_k=\frac{1}{(n-k)!}p^{n-k}$ for $k=1,\dots,n$, then the only non-zero commutation
relations are: $[X_0, X_{k}]=X_{k+1}$ for $k=1,\dots,n-1$. 
\item [{\bf E7}] The $(n+2)$ elements $pq, q, 1, p, \ldots, p^{n-1}$
 span a non-nilpotent solvable Lie subalgebra whose derived algebra is isomorphic to
 ${\cal L}_n$.  If we set $h=pq, X_0=-q$ and
 $X_k=\frac{1}{(n-k)!}p^{n-k}$ for $k=1,\dots,n$, then the only non-zero commutation
relations are: $[h,X_0]=X_0, [h, X_k]=-(n-k) X_k$ and $ [X_0, X_{k}]=X_{k+1}$
 for $k=1,\dots,n-1$. We denote this Lie algebra by $\tilde {\cal L}_n$. 
\item [{\bf E8}] The $(n+1)$ elements $pq, p^{i_1}, \ldots, p^{i_n}$, where $i_1,\dots, i_n$ are  distinct positive integers not all zero, 
 span a non-nilpotent solvable Lie subalgebra whose derived algebra is
  $n-$dimensional and abelian.  
We denote this Lie algebra by $\rr(i_1,\dots,i_n)$. 
It is clear that $\rr(i_{\sigma(1)},\dots,i_{\sigma(n)})\cong \rr(i_1,\dots,i_n)$ for any permutation $\sigma\in S_n$, that $\rr(i_1,\dots,i_n)$ has a non-trivial centre iff one of the indices is zero and that $\rr(0,i_2,\dots,i_n)\cong {\mathfrak r}_{n-1}(i_2,\dots,i_n)\times\CC$. Note also that $\rr (mi_1,\dots mi_n)\cong \rr (i_1,\dots i_n)$ if $m\in\NN^*$. If we set $h=pq, X_k=p^{i_k}$ for $k=1,\dots,n$ the only non-zero commutation relations are: $[h, X_k]=-i_k X_k$. 
\end{enumerate}

\section{Classification of finite-dimensional Lie algebras which can be realised in the Weyl algebra}
\renewcommand{\theequation}{4.\arabic{equation}}   
\setcounter{equation}{0}

In this section we obtain all finite-dimensional Lie algebras that can
be realised as subalgebras of $A_1$.
The only such non-solvable Lie algebras  are isomorphic to either  $\s$, $\s\times\CC$ or $\s\ltimes {\cal H}_3$.
This is basically a consequence of the fact that $A_1$ satisfies {\bf P3},  the 
``commutative centraliser condition''.
We then show that a non-abelian nilpotent Lie algebra which can be realised as a subalgebra of $A_1$ is
isomorphic to an ${\cal L}_n$  (cf {\bf E6})
and this is consequence of the properties  {\bf P3} and {\bf P3}${}^\prime$.
Finally, we
show that a solvable non-nilpotent Lie algebra which can be realised
as a subalgebra of $A_1$ is isomorphic to either an  $\tilde {\cal L}_n$ (cf {\bf E7}) or to  an
$\rr(i_1,\dots,i_n)$  (cf {\bf E8}).
This result is more difficult to prove and follows from special 
properties of the spectrum of semi-simple elements of $A_1$.

\subsection{Non-solvable Lie algebras}

\begin{proposition}
\label{simple}
(Theorem $3$ of \cite{italienii}) If $\mathfrak{g}$ is a semi-simple complex Lie algebra of rank $>1$, then $A_1^{\mathfrak{g}}=\emptyset$.
\end{proposition}
\textbf{Proof}: Suppose for contradiction that there exists $f\in
 A_1^{\mathfrak{g}}$.
Let  $\mathfrak{h}\subset\mathfrak{g}$ be a Cartan subalgebra and let
$$\mathfrak{g}=\mathfrak{n}_- \oplus \mathfrak{h} \oplus \mathfrak{n}_+$$
be the corresponding triangular decomposition for some choice of simple
roots.
 Let $x\in\mathfrak{n}_+$
be a non-zero highest root vector. 
 The commutant of $x$
in  $\mathfrak{g}$, $Z_{\mathfrak{g}}(x)$, contains $\mathfrak{n}_+$
which is not abelian 
 since  rank$(
{\mathfrak{g}})>1$. 
Hence 
 $C(f(x))$ contains $f(\mathfrak{n}_+)$ which  is not abelian and so by the \textit{ccc}, $f(x)\in \CC$ and $x\in Z_\g$. But $\g$ is semi-simple and hence $x=0$  which is a contradiction. QED

\begin{proposition}
\label{reductive}
Let $\g=\g_1\times\mathfrak z$ be a reductive complex Lie algebra where
$\g_1$ is  semi-simple  and  $\mathfrak z$ is the centre.
If $A_1^{\mathfrak{g}}\ne\emptyset$, then
\begin{enumerate}
\item $\g_1\cong \s$.
\item For any $f\in A_1^{\mathfrak{g}}$, $f({\mathfrak z})\subseteq\CC$. 
\end{enumerate}
\end{proposition}
\textbf{Proof}: By  Proposition \ref{simple}  one has $\g_1\cong\s$.
 Let $z$ be an element of $\mathfrak z$. Then $C(f(z))$ contains $f(\g_1)$ which is not abelian and hence, by the \textit{ccc}, $f(z)$ is a scalar. QED

\begin{remark}
A noncommutative algebra $A$ is said to satisfy the commutative centraliser condition
if the centraliser of any element not in the centre $Z_A$ 
 is a commutative subalgebra. Proposition \ref{reductive} then remains true if we replace $A_1$ by $A$ and $\CC$ by $Z_A$.
One example of such an algebra is the universal enveloping algebra ${\cal U}(\s)$ (for other examples see \cite{Bavula}).
\end{remark}

\begin{remark}
\label{super}
Using the classification theorem of Kac,  this result implies that a classical complex simple Lie superalgebra contained in $A_1$ is isomorphic to ${\mathfrak {osp}}(1\vert2)$.
\end{remark}

By the theorem of Levi-Malcev, any finite-dimensional non-solvable Lie algebra $\g$ is isomorphic 
to a semi-direct product of  its radical $\rr$ and a semi-simple subalgebra $\mathfrak s$.  Suppose $\g$ is realisable in $A_1$. Then  $\mathfrak s$ is isomorphic to $\s$ by Proposition \ref{simple} and
 by analysing the action of this  $\s$ on $\rr$ we will show that there are in fact only three possibilities for $\rr$.

\begin{definition}
\label{triplet}
(i) Three non-zero elements $X,Y,H$ of $A_1$ are called an \textit{$\s$ triplet} if they satisfy the relations:
\beqa
[H,X]=2X,\quad
[H,Y]=-2Y,\quad
[X,Y]=H.
\eeqa
(ii) An element $v\in A_1$ is of weight $\lambda\in\CC$ if
$[H,v]=\lambda v$.\\
(iii) 
The set of elements  of weight $\lambda$ in a linear 
subspace $E\subseteq A_1$ will be denoted by
$E_\lambda$. (Note that  $E_0$ is abelian  by the {\it ccc}.)
\end{definition}

\begin{lemma}
\label{lema-marcus}
Let $X, Y, H$ be an $\s$ triplet and
let $\mathfrak l\subseteq A_1$ be a Lie subalgebra stable under 
${\mathrm {ad}}(X),{\mathrm {ad}}(Y)$ and ${\mathrm {ad}}(H)$. Suppose there
exists $v\ne 0$
 in $\mathfrak l$ of weight  $\lambda\in\CC^*$ such that 
$[X,v]=0$. Then $[v, [Y,v]]\ne 0$ and is of weight $2\lambda -2$. 
\end{lemma}
\textbf{Proof:} We suppose for contradiction that $[v, [Y,v]]=0$. Then
$C(v)$ contains $X$ and $[Y,v]$. But
$[X,[Y,v]]=[H,v]+[Y,[X,v]]=\lambda v\ne 0$ and therefore $C(v)$
is non-abelian. By the \textit{ccc}, $v\in \CC$ and so $[H,v]=0$ which is a contradiction.
 Hence, $[v, [Y,v]]\ne 0$ and is
obviously of
weight $2\lambda -2$.
QED

\begin{proposition}
\label{oarecare}
Let $X, Y, H$ be an $\s$ triplet and
let ${\mathfrak l}\subseteq A_1$ be a finite-dimensional  Lie
subalgebra stable under  ${\mathrm {ad}}(X),{\mathrm {ad}}(Y)$ and ${\mathrm {ad}}(H)$. 
Then ${\mathfrak l}_\lambda=\{ 0 \}$ for $|\lambda|>2$.
\end{proposition}
\textbf{Proof:} Recall first that since $\mathfrak l$ is a \textit{finite}-dimensional representation of $\s$, ${\mathfrak l}_\lambda =\{ 0 \}$
iff ${\mathfrak l}_{-\lambda} =\{ 0 \}$.
Let $\lambda_{max}$ be the largest eigenvalue of $\mathrm{ad}(H)$
restricted to ${\mathfrak l}$ and suppose  that $\lambda_{max}>2$. 
Then $\lambda_{max}\ge2\lambda_{max}-2$
by Lemma \ref{lema-marcus} which is a contradiction. QED

\begin{proposition}
\label{solvable}
Let $X, Y, H$ be an $\s$ triplet and
let ${\mathfrak r}\subseteq A_1$ be a finite-dimensional  solvable Lie subalgebra stable under  ${\mathrm {ad}}(X),{\mathrm {ad}}(Y)$ and ${\mathrm {ad}}(H)$. 
Then ${\mathfrak r}_\lambda=\{ 0 \}$ for $\vert\lambda\vert>1$ and $\rr_0\subseteq\CC$.
\end{proposition}
\textbf{Proof:} Recall that if $\g$ is a Lie algebra the upper central series $(\g^{(i)})_{i\in\NN}$ is
defined 
inductively by:
$\g^{(0)}=\g, \ \g^{(i+1)}=[\g^{(i)},\g^{(i)}].$
 Suppose for contradiction that ${\mathfrak r}_2\ne \{ 0 \}$.
Then, by  Lemma \ref{lema-marcus}, ${\mathfrak r}^{(i)}_2\ne \{ 0 \}$
for all positive $i$. But ${\mathfrak r}$ is solvable so that by
definition ${\mathfrak r}^{(m)}=\{ 0 \}$ for some $m\in \NN$ and this is a
contradiction. Now, since $\rr={\mathfrak r}_0\oplus{\mathfrak r}_{-1}\oplus{\mathfrak r}_1$ it follows that $[X, \rr_0]=[Y,\rr_0]=\{ 0 \}$ and therefore, by the {\it ccc}, that $\rr_0 \subseteq \CC$.
QED

\begin{proposition}
\label{finala}
Let $X, Y, H$ be an $\s$ triplet and
let ${\mathfrak r}\subseteq A_1$ be a finite-dimensional  solvable Lie
subalgebra stable under ${\mathrm ad}(X),{\mathrm ad}(Y)$ and ${\mathrm
ad}(H)$. 
Then ${\mathfrak r}$ is isomorphic to either $\{ 0 \}$, $\CC$ or ${\cal H}_3$
 (the three-dimensional Heisenberg algebra).
\end{proposition}
\textbf{Proof:} 
We have ${\mathfrak r}={\mathfrak r}_0\oplus{\mathfrak r}_{-1}\oplus{\mathfrak r}_1$ (by Proposition \ref{solvable}), $ {\mathfrak r}_0\subseteq \CC$ (by Proposition \ref{solvable}) and $[{\mathfrak r}_{-1},{\mathfrak
r}_1]\subseteq{\mathfrak r}_0$ since ${\mathrm
{ad}}(H)$ is a derivation.

Suppose first that $\mathrm{dim}\ {\mathfrak r}_{1}\ge 2$ and let
$v\in{\mathfrak r}_1\setminus \{0\}$. 
The kernel, $Z_{{\mathfrak r}_{-1}}(v)$, of the linear map 
$\mathrm{ad}(v):{\mathfrak r}_{-1}\to{\mathfrak r}_{0}$ is of dimension $\ge 1$ and therefore contains a non-zero vector $w$. Hence $C(v)$ contains $X$ and $w$. But $[X,w]\ne 0$ since $\mathrm{ad}(X) : {\mathfrak r}_{-1}\to {\mathfrak r}_{1}$ is an isomorphism and so $C(v)$ is not abelian. 
By the \textit{ccc}, $v$ is a scalar which is a contradiction since $[H,v]=v$.
Therefore $\mathrm{dim}\ {\mathfrak r}_{-1}=\mathrm{dim}\ {\mathfrak r}_{1}\le 1$.

If $\mathrm{dim}\ {\mathfrak r}_{1}=1$, then $ {\mathfrak r}$ is not abelian by 
Lemma \ref{lema-marcus} and so ${\mathfrak
r}_0=\CC$ by Proposition \ref{solvable}. 
Hence
${\mathfrak r}_0$ is the centre of ${\mathfrak r}$ and it is now
obvious that ${\mathfrak r}$ is isomorphic to the three dimensional Heisenberg algebra.
 If $\mathrm{dim}\ {\mathfrak r}_{1}=0$, then  ${\mathfrak r}_{0}=\{0\}$ or $\CC$  by Proposition \ref{reductive}. QED

\bigskip

We can now conclude  this subsection with the following theorem:
\begin{theorem}
\label{mareata}
Let $\g$ be a finite-dimensional non-solvable Lie algebra. Then $A_1^\g\ne\emptyset$ iff $\g$ is isomorphic to one of the following:
\begin{enumerate}
\item $\s$,
\item $\s \times \CC$,
\item $\s\ltimes\mathcal H_3 $.
\end{enumerate}
\end{theorem}
\textbf{Proof:} 
Immediate from 
 the Levi-Malcev theorem, Propositions \ref{simple} and \ref{finala}.
Note that by the {\it ccc}, $\s\times\mathcal H_3 $ cannot be realised as a Lie subalgebra of $A_1$. QED

\medskip

\begin{corollary}
\label{KM}
Let $\g$ be a complex Lie algebra which contains a finite-dimensional
non-solvable subalgebra not isomorphic to $\s$, $\s\times \CC$ or
$\s\ltimes\mathcal H_3 $. Then $A_1^{\mathfrak{g}}=\emptyset$.
\end{corollary}

\begin{corollary}
\label{afine}
Let $\tilde \g$ be the affine Kac-Moody algebra associated to the
simple complex Lie algebra $\g$. Then $A_1^{\tilde {\mathfrak{g}}}=\emptyset$.
\end{corollary}
\textbf{Proof:} Since $\tilde \g$ contains a
reductive Lie algebra isomorphic to $\g\times \CC^2$ (generated by $\g$, the derivation and the
central element), the result follows immediately from Corollary \ref{KM}. QED

\subsection{Nilpotent non-abelian Lie algebras}

\begin{theorem}
\label{ln}
Let $\n\subset A_1$ be a nilpotent, non-abelian Lie
subalgebra of dimension $n$. Then $\n\cong {\mathcal L}_{n-1}$.
\end{theorem}
{\bf Proof:} Let
$(\n_{(i)})_{i\in\NN}$ defined by
$\n_{(0)}=\n$ and $\n_{(i+1)}=[\n,\n_{(i)}]$ be the lower central series. Let  $k\ne 0$ be the unique
positive integer
such that $\n_{(k)}\ne\{ 0 \}$ and $\n_{(k+1)}=\{ 0 \}$. The theorem will essentially be a consequence of the following lemma:

\begin{enumerate}
\item[]
\begin{lemma}
\label{nilpotentele}
\begin{enumerate} 
\item $Z_\n=\n_{(k)}=\CC$.
\item There exist $P,Q\in \n$ such that $[P,Q]=1$ and $\n=<P>\oplus\ 
Z_\n (Q)$. 
\item $Z_\n (Q)$ is abelian and $\n_{(1)}\subset Z_\n (Q)$.
\item ${\mathrm{dim}}\ \n_{(i)}=n-i-1$
for $1\le i\le k$.
\item $k=n-2$.
\end{enumerate}
\end{lemma}
\textbf{Proof:} $(a)$: Let $z\in Z_\n$. Then $C(z)$ contains
$\n$ which is not commutative. By the \textit{ccc}, $z$ is a
scalar and therefore $Z_\n\subseteq\CC$. But
 $Z_n$ contains $\n_{(k)}$ and so $Z_\n = \n_{(k)}=\CC$.

$(b)$: Since $\n_{(k)}=[\n,\n_{(k-1)}]=\CC$, there exist $P\in\n,\
Q\in\n_{(k-1)}$ such that $[P,Q]=1$. 
It is clear that  $<P>\cap \ Z_\n
(Q)=\{ 0 \}$ since $P$ and $Q$ do not commute. We now show that $\n=<P>+\ 
Z_\n (Q)$.
Let $v\in\n$. There exists $\lambda\in\CC$ such that $[v,Q]=\lambda$ since
$[\n,\n_{(k-1)}]=\CC$. Then $v-\lambda P\in Z_\n(Q)$ and $v=\lambda
P+(v-\lambda P)\in\  <P>+\ Z_\n (Q)$.

$(c)$: By the \textit{ccc}, $Z_\n (Q)$ is abelian. Let $z\in Z_\n(Q)$. Then 
$$ [Q,[P,z]]=[[Q,P],z]+[P,[Q,z]]=0$$
because $[Q,P]=-1$ and $[Q,z]=0$. Hence $[P, Z_\n (Q)]\subseteq Z_\n(Q)$ and since $[Z_\n (Q), Z_\n (Q)]=\{ 0 \}$, 
this means that $\n_{(1)}=[\n,\n]$ is contained in $Z_\n (Q)$. 

$(d)$:   Since $\n=<P>\oplus\ 
Z_\n (Q)$ and since $Z_\n (Q)$ is abelian, we have
\beqa
\label{adjoainta}
\n_{(i)}=\mathrm{ad}^i(P)(Z_\n(Q)) \ \forall i\in\NN^*.
\eeqa
\noi
If $z\in Z_\n (Q)$ is such that $\mathrm{ad}(P)(z)=0$ then $z\in
 C (P)\cap C (Q)$. By property {\bf P3}${}^{\prime}$,
 $C (P)\cap C (Q)=\CC$ since $P$ and $ Q$ do not commute.
Hence $z\in\CC$ and 
\beqa
\mathrm{Ker}\ \mathrm{ad}(P) \vert _{\n}=<P>\oplus\ \CC,\quad
\mathrm{Ker}\ \mathrm{ad}(P) \vert _{Z_\n(Q)}=\CC.
\eeqa
This implies $(d)$. 

$(e)$: It follows immediately from $(d)$ that $k=n-2$. QED
\end{enumerate}

We will now prove that $\n\cong {\cal L}_{n-1}$.
Since $\n_{(k-1)}=\mathrm{ad}^{k-1}(P)(Z_\n(Q))$ and
since $Q\in \n_{(k-1)}$, there exists $w\in Z_\n(Q)$ such that
$\mathrm{ad}^{k-1}(P)(w)=Q$. It is clear that $X_1=w,
X_2=\mathrm{ad}(P)(w),X_3=\mathrm{ad}^2(P)(w),\ldots,X_{k+1}=\mathrm{ad}^{k}(P)(w)$ 
are linearly independent (since
$\mathrm{ad}^{k}(P)(w)\ne 0$ and 
$\mathrm{ad}^{k+1}(P)(w)=0$) and by a dimension count, these vectors
form a basis of $Z_\n(Q)$. The only non-zero commutation relations of $\n$ in
the basis $X_0=P, X_1, \ldots, X_{n-1}$ are 
$$[X_0,X_i]=X_{i+1} \mbox{ for }  1\le i \le n-2, $$
which are the standard commutation relations for ${\cal
L}_{n-1}$ (cf {\bf E6}). QED

\subsection{Solvable non-nilpotent Lie algebras}

\begin{lemma}
\label{solvable1}
Let ${\mathfrak g}\subset A_1$ be a finite-dimensional 
 Lie subalgebra and let $\g'$ be its derived algebra. Suppose there exists  $h\in \g$ be such that
$\mathrm{ad}(h)\vert_{{\mathfrak g}}$ is not nilpotent.  Let $0,\lambda_1,\dots,\lambda_k$ be the
 distinct eigenvalues of  $\mathrm{ad}(h)\vert_{{\mathfrak g}}$ and let $E_{0},
 E_{\lambda_1},\dots,E_{\lambda_k}$ be the corresponding eigenspaces.
\begin{enumerate}
\item $\g=E_{0}\oplus E_{\lambda_1}\oplus \ldots \oplus E_{\lambda_k}$
\item $\g'=(\g'\cap E_{0})\oplus E_{\lambda_1} \oplus\ldots \oplus E_{\lambda_k}$
\item $E_0=<h>$ or $E_0=<h>\oplus\ \CC$.
\end{enumerate}
\end{lemma}
\textbf{Proof:}
By the theory of endomorphisms,
$$ {\mathfrak g}=\mathrm{Ker}\  \mathrm{ad}^{m_0}(h)\vert_{{\mathfrak g}}\oplus\mathrm{Ker} (\mathrm{ad}({h})\vert_{{\mathfrak g}} -
\lambda_1)^{m_1}\oplus \dots \oplus \mathrm{Ker} (\mathrm{ad}({h})\vert_{{\mathfrak g}} -
\lambda_k)^{m_k},$$ 
\noi where the characteristic polynomial of
$\mathrm{ad}({h})\vert_{{\mathfrak g}}$ is
$x^{m_0}(x-\lambda_1)^{m_1}\ldots(x-\lambda_k)^{m_k}$.
By Proposition $6.5$ of \cite{Dixmier},
$$\mathrm{Ker} (\mathrm{ad}({h}) -
\lambda_j)^{m_j}=\mathrm{Ker} (\mathrm{ad}({h}) -
\lambda_j).$$
Since $\mathrm{ad}({h})$ has a non-zero
eigenvalue, $C(h)\ne D(h)$ (cf subsection $2.2$) 
and therefore, by Theorem \ref{th-Dixmier}, $C(h)=N(h)$ so that 
$$\mathrm{Ker} \ \mathrm{ad}^{m_0}(h)=\mathrm{Ker} \ \mathrm{ad}({h}).$$
This proves part $1$.
To prove part $2$ it is sufficient to note that
$E_{\lambda_i}\subseteq \g'$ since 
$v_i=\frac{1}{\lambda_i} [h, v_i]$ for all $v_i\in E_{\lambda_i}$. 

To prove part $3$, let $h'\in E_0$. Then $[h',E_{\lambda_i}]\subseteq
E_{\lambda_i}$ since $[h,h']=0$. Thus there exist $\alpha\in\CC$,
$v\in E_{\lambda_1}\setminus \{ 0 \}$ such that $[h',v]=\alpha v$. Hence
$[h'-\frac{\alpha}{\lambda_1}h, v]=0$ which means that $h$ and $v$
commute with $h'-\frac{\alpha}{\lambda_1}h$. But $[h,v]\ne 0$ and so by the {\it
ccc},  $h'-\frac{\alpha}{\lambda_1}h\in\CC$. QED

\begin{remark}
\label{simplificare}
Taking $\g$ semi-simple, this lemma provides an alternative proof of
Proposition \ref{simple}.
However the proof we gave for Proposition \ref{simple} works for any
algebra satisfying the {\it ccc}, whereas the proof of the theorem above
depends on the existence of the Dixmier partition (Theorem \ref{th-Dixmier}) and other special properties not in general available for an algebra satisfying the {\it ccc}.
\end{remark}

If $\g\subset A_1$ is a finite-dimensional solvable 
non-nilpotent Lie subalgebra then, by Theorem \ref{ln},  
the derived algebra $\g'$ is  isomorphic to either an ${\cal
L}_n$ or  is abelian. We first treat the case where $\g'\cong {\cal
L}_n$.

\begin{theorem}
\label{solvable22}
Let ${\mathfrak r}\subset A_1$ be a finite-dimensional solvable 
non-nilpotent Lie subalgebra whose derived algebra ${\mathfrak
r}'$ is isomorphic to ${\cal L}_n$. 
Then  ${\mathfrak r}$ is
isomorphic to $\tilde {\cal L}_n$.
\end{theorem}
{\bf Proof:} 
Since ${\mathfrak r}\subset A_1$ is a finite-dimensional  non-nilpotent Lie algebra, there
exists $h\in {\mathfrak r}$  such that $\mathrm{ad}(h)\vert_{{\mathfrak r}}$ is not nilpotent by 
Engel's theorem.
By Lemma \ref{solvable1}, $\mathrm{ad}(h)\vert_{{\mathfrak r}'}$ is
diagonalisable and therefore
by Theorem $1$ of \cite{GOZE} there exists a basis $X_0,\ldots, X_n$ of
${\mathfrak r}'\cong {\cal L}_n$ of eigenvectors of $\mathrm{ad}(h)$ such that the only
non-zero commutation relations are
\beqa
\label{baza}
[X_0, X_i]&=&X_{i+1},\ i=1,..,n-1\nonumber \\
\left[ h, X_j\right]&=&\alpha_j X_j,\  j=0,..,n.
\eeqa
\noi
Since $\CC\subseteq\rr'$ by Theorem \ref{ln},  we must have 
$X_n\in\CC$,  $[h,X_n]=0$ and $\alpha_n=0$. 
But $\mathrm{ad}(h)\vert_{{\mathfrak r}'}$ is a (non-zero) derivation so its eigenvalues 
 satisfy 
\beqa
\label{cond1}
\alpha_{i+1}=\alpha_0+\alpha_i, \ i=1,..., n-1.
\eeqa
\noi
From this it follows that the eigenvalues of $\mathrm{ad}({h})|_{\rr'}$
are $\alpha_0, (1-n)\alpha_0,(2-n)\alpha_0,\dots, 0$ and hence
$\alpha_0\ne 0$ and  the
eigenspaces are of dimension $1$. The only non-zero commutation relations of
$\rr$ are then
\beqa
\label{baza-fin}
[X_0, X_i]&=&X_{i+1},\ i=1,..,n-1 \nonumber \\
\left[ \frac{1}{\alpha_0} h, X_0\right]&=& X_0 \nonumber \\
\left[  \frac{1}{\alpha_0} h, X_i\right]&=&- (n-i)X_i \ i=1,..,n-1
\eeqa
 and this shows that $\rr$ is isomorphic to $\tilde {\cal{L}}_n$ (cf {\bf E7}). QED 

\bigskip

It now remains to treat the case where $\g\subset A_1$ is a finite-dimensional solvable non-nilpotent Lie subalgebra whose derived algebra is abelian. 
We will need the following lemma (also see Theorem $3.2$ of \cite{joseph-jmp}):

\begin{lemma}
\label{joseph}
Let $h, X_1, X_2\in A_1\setminus\{ 0 \}$ and
 $(\lambda_1,\lambda_2)\in \ZZ^2\setminus \{(0, 0)\}$
be such that 
\beqa
\label{cond-joseph}
[h, X_1]=\lambda_1 X_1,\ [h, X_2]=\lambda_2 X_2,\mbox{ and }[X_1,
X_2]=0.
\eeqa
\noi Then $\lambda_1 \lambda_2 > 0$ and
there exists $a\in\CC^*$ such that
$X_1^{|\lambda_2|}=a X_2^{|\lambda_1|}$.
\end{lemma}
{\bf Proof}: Since $[X_1, X_2]=0$, by Theorem $3.1$ of \cite{joseph-jmp} there exist
$m,n\in\NN^*$ and
$\alpha_{ij}\in\CC$ not all equal to zero such that
\beqa
\label{lema-joseph}
\sum_{i=0}^m \sum_{j=0}^n \alpha_{ij} X_1^i X_2^j=0
\eeqa
This can be rewritten as
\beqa
\sum_{u\in U} \left( \sum_{(i,j): i \lambda_1 + j \lambda_2=u} \alpha_{ij} X_1^i X_2^j\right)=0
\eeqa
\noi
where
$ U=\{ i \lambda_1 + j \lambda_2 \in \ZZ : \ 0\le i \le m, 0\le j \le n\}$.
Since eigenvectors corresponding to distinct eigenvalues of $h$ are linear
independent we deduce that for all $u\in U$, 
$$\sum_{(i,j): i \lambda_1 + j \lambda_2=u} \alpha_{ij} X_1^i X_2^j=0.$$
Choose $(i_0,j_0)\in \ZZ^2$ such that $\alpha_{i_0j_0}\ne 0$. We
set $u_0=i_0 \lambda_1 + j_0 \lambda_2$
and 
$S=\{(i,j)\in \ZZ^2:\ i \lambda_1 + j \lambda_2=u_0\}$. 
In $\RR^2$ the solutions $(x,y)$ of the equation $x \lambda_1 +
y\lambda_2=u_0$ define an affine line and $S$ is a discrete subset of
this line.
It is then easy to see that 
there exist $(i',j')\in \ZZ^2$ and $(i_m,j_m)\in S$ such that $j'\ge 0$,
\beqa
\label{aj-simplu}
i'\lambda_1+j'\lambda_2=0
\eeqa
\noi
and every
element $(i,j)$ of $S$ is of the form $(i,j)=(i_m, j_m)+k_{ij}(i',j')$ for
some $k_{ij}\in\NN$. 
One can then write (in the field of fractions of $A_1$) 
$$\sum_{(i,j): i \lambda_1 + j \lambda_2=u_0} \alpha_{ij} X_1^i
X_2^j=X_1^{i_m}X_2^{j_m} \sum_{(i,j): i \lambda_1 + j \lambda_2=u_0}\alpha_{ij} (X_1^{i'}X_2^{j'})^{k_{ij}}=0.$$
\noi
We thus have a polynomial $P$ in the variable
$X_1^{i'}X_2^{j'}$ and factorising (in the field of fractions of
$A_1$) we deduce that there exists $c\in\CC^*$ such that
$X_1^{i'}X_2^{j'}=c$. Hence
$X_1^{i'\lambda_1}X_2^{j'\lambda_1}=c^{\lambda_1}$,  $(X_1^{-\lambda_2}X_2^{\lambda_1})^{j'}=c^{\lambda_1}$ (by
(\ref{aj-simplu})) and therefore $X_1^{-\lambda_2}X_2^{\lambda_1}\in\CC^*$. 

If $\lambda_1 \lambda_2\le 0$ this means that
there exists $b\in\CC^*$ such
that $ X_1^{|\lambda_2|}X_2^{|\lambda_1|}=b$ which is clearly
impossible since in $A_1$ the only invertible elements are the
scalars (see page $210$ of \cite{Dixmier}). Hence $\lambda_1\lambda_2>0$ and there exists $a\in\CC$ such that $X_1^{|\lambda_2|}=a X_2^{|\lambda_1|}$. QED

\medskip

We can now show that  a finite-dimensional solvable
non-nilpotent Lie subalgebra of $A_1$ whose derived algebra is abelian is isomorphic to an 
 $\rr(i_1,\dots i_n)$ (cf  {\bf E8}).

\begin{theorem}
\label{solvable2}
Let ${\mathfrak r}\subset A_1$ be a finite-dimensional solvable 
non-nilpotent Lie subalgebra whose derived algebra ${\mathfrak r}'$ is abelian. 
Then there exist distinct positive integers $i_1,\ldots,i_n$ not all zero such that ${\mathfrak r}$ is
isomorphic to $\rr(i_1,\dots i_n)$.
\end{theorem}
{\bf Proof:}  Let $h\in {\mathfrak r}$ be such that
$\mathrm{ad}(h)\vert_{{\mathfrak r}}$ is not nilpotent, let $0,\lambda_1,\dots,\lambda_k$ be its
 distinct eigenvalues and let $E_{0},
 E_{\lambda_1},\dots,E_{\lambda_k}$ be the corresponding eigenspaces.
By Lemma \ref{solvable1},
$$\rr=E_0\oplus E_{\lambda_1}\oplus E_{\lambda_2} \oplus \dots \oplus E_{\lambda_k}$$
and since $ E_{\lambda_1}\oplus E_{\lambda_2} \oplus \dots
E_{\lambda_k}\subseteq \rr'$ is abelian we deduce that
$$\rr'=E_{\lambda_1}\oplus E_{\lambda_2} \oplus \dots \oplus E_{\lambda_k}.$$

It is clear that if the centre of $\rr$ is $\CC$, then $\rr\cong
 \tilde \rr\times\CC$ where $\tilde \rr$ satisfies the hypothesis of
the theorem and has trivial centre.
Therefore (see {\bf E8}) to prove the theorem  it is enough to consider the case
where $\rr$ has trivial centre.
 
First, note that by Lemma \ref{solvable1}($3$) we have  $E_0=<h>$. Next,
since $h\in\Delta_3 \cup \Delta_4$ (by Remark \ref{simplificari}),
there exists $\rho\in\CC^*$ such that eigenvalues of $\mathrm{ad}(h)$
are integer multiples of $\rho$ (Corollary $9.3$ of \cite{Dixmier} if
$h\in\Delta_3$ and Theorem $1.3$ of \cite{joseph3} if
$h\in\Delta_4$). Hence $\frac{1}{\rho}\mathrm{ad}(h)\vert_{\rr}$ has integer eigenvalues
 and, since the $E_{\lambda_i}$ are abelian,  these integers are all of the same sign
 and the corresponding eigenspaces are of dimension $1$ by Lemma \ref{joseph}.
It is now clear that $\rr$ is isomorphic to 
$\rr(|\frac{\lambda_1}{\rho}|,\dots , |\frac{\lambda_k}{\rho}|)$. QED

\subsection{Finite-dimensional Lie subalgebras of Der$(A_1)$}

In this subsection we find all finite-dimensional  Lie
algebras that can be realised as subalgebras of Der$(A_1)$ and 
we give some examples of Lie subgroups of $\A (A_1)$ which exponentiate them.

\begin{theorem}
\label{ptintro}
Let $\g\subseteq \mathrm{Der}(A_1)$ be a finite-dimensional non-abelian Lie subalgebra. Then $\g$ is  isomorphic to one of the following: 
\beqa
&& 1) \ \s, \nonumber \\
&& 2) \ \s \ltimes \CC^2, \nonumber \\
&& 3) \ {\cal L}_n\ (n\ge 2),\nonumber \\
&& 4) \ \tilde {\cal L}_n/\CC \ (n\ge 2), \nonumber \\
&& 5) \ \rr (i_1,\dots ,i_n) \ \ (0<i_1<\dots <i_n). \nonumber
\eeqa
\end{theorem}
{\bf Proof:} Consider the commutative diagram:
$$\xymatrix{A_1\ar[r]^{\pi}&A_1/{\mathbb C}\ar[r]^{{\mathrm{ad}}}
&{\rm Der}(A_1)\\
\pi^{-1}({\rm ad}^{-1}(\mathfrak g))\ar@{^{(}->}[u]\ar[r]^{\pi}
&{\rm ad}^{-1}(\mathfrak g)\ar@{^{(}->}[u]\ar[r]^{\rm ad}
&\mathfrak g\ar@{^{(}->}[u]
}$$ 

In this diagram, $\pi: A_1\to A_1/\CC$ is a surjective Lie algebra
homomorphism and \\
$\mathrm{ad}:A_1/\CC\to \mathrm{Der}(A_1)$ is a Lie algebra isomorphism (see page
$210$ of \cite{Dixmier}). 

If $\g$ is not abelian then $\pi^{-1}(\mathrm{ad}^{-1} (\g)) $ is a finite-dimensional
non-abelian subalgebra of $A_1$ containing $\CC$ and $\g\cong \pi^{-1}(\mathrm{ad}^{-1}
(\g)) / \CC$.
By the results of the previous subsections $\pi^{-1}(\mathrm{ad}^{-1}
(\g)) $ is isomorphic to one of the following: $\s\times\CC$, $\s
\ltimes {\cal H}_3$, an ${\cal L}_n$, an $\tilde {\cal {L}}_n$ or an ${\mathfrak
r} (0, i_2,\dots,i_n)$.
Taking the quotient by $\CC$ this shows that $\g$ is isomorphic to one
of: $\s$, $\s\ltimes \CC^2$, an ${\cal {L}}_n / \CC$, an
$\tilde {\cal {L}}_n/ \CC$ or an $\rr(i_2,\dots,i_n)$ (cf {\bf E8}). Since ${\cal {L}}_n / \CC\cong {\cal {L}}_{n-1}$ 
 this proves the theorem. QED

\begin{remark}
The derived Lie algebra of $\tilde {\cal {L}}_n/ \CC$ is isomorphic to ${\cal L}_{n-1}$ but $\tilde {\cal {L}}_n/ \CC$ is isomorphic to neither $\tilde {\cal {L}}_{n-1}$ nor  ${\cal L}_{n}$, both of whose derived algebras are also isomorphic to ${\cal L}_{n-1}$. A basis for $\tilde {\cal {L}}_n/ \CC$ is $\bar h, \bar X_0, \dots, \bar X_{n-1}$ (see {\bf E7} for notations) and the only non-zero commutation relations are  $[\bar h,\bar X_0]=\bar X_0$, $[\bar h, \bar X_k]=-(n-k)\bar X_k$ for $k=1,\dots,n-1$ and $[\bar X_0, \bar X_k]=\bar X_{k+1}$ for $k=1,\dots,n-2$. In particular, $\tilde {\cal {L}}_n/ \CC$ is  not nilpotent and its centre is trivial.
\end{remark}

This theorem implies that if $G$ is a finite-dimensional,
connected and  in some sense Lie subgroup
of $\A (A_1)$, 
then $G$ is either abelian or a discrete quotient of
$SL(2,\CC)$, 
$SL(2,\CC)\ltimes \CC^2$ or of the simply-connected Lie groups $L_n$, $\tilde L_n/\CC$ and  $R (i_1, \ldots, i_n)$ corresponding respectively to the Lie algebras ${\cal L}_n$, $\tilde {\cal L}_n/\CC$ and  $\rr (i_1,\dots, i_n)$.
We now give explicit constructions of $R (i_1, \ldots, i_n)$, $\tilde L_n$ and $L_n$ and show that the groups
 $SL(2,\CC)$, 
$SL(2,\CC)\ltimes \CC^2$, $R (i_1, \ldots, i_n)/\ZZ$, $\tilde L_n/(\CC\times\ZZ)$ and $L_n$ can be holomorphically embedded
 in $\A(A_1)$.

\begin{enumerate}
\item [{\bf E9}] Define 
$\hat \alpha_1 : SL(2,\CC)\to \A (A_1)$ by 
$$\hat
\alpha_1 ( \begin{pmatrix} a_1 & a_2 \\ a_3 & a_4 \end{pmatrix} )
(p)=a_2 q +a_4 p,\ \  \hat
\alpha_1 ( \begin{pmatrix} a_1 & a_2 \\ a_3 & a_4 \end{pmatrix} )
(q)=a_1 q +a_3 p$$
 and $\hat \alpha_2 : \CC^2\to \A (A_1)$   
 by  
$$\hat
\alpha_2 ( \begin{pmatrix} b_1 \\ b_2 \end{pmatrix} )
(p)=p-b_1,\ \  \hat
\alpha_2 ( \begin{pmatrix} b_1 \\ b_2 \end{pmatrix} )
(p)=q+b_2.$$   
Then one checks that  $\hat \alpha_1$ and $\hat \alpha_2$
are injective group homomorphisms such that
$$ \hat \alpha_2 ( \begin{pmatrix} a_1 & a_2 \\ a_3 & a_4
\end{pmatrix} \begin{pmatrix} b_1 \\ b_2 \end{pmatrix} ) = \hat
\alpha_1 ( \begin{pmatrix} a_1 & a_2 \\ a_3 & a_4 \end{pmatrix} ) \ \hat
\alpha_2 ( \begin{pmatrix} b_1 \\ b_2 \end{pmatrix} ) \  \hat
\alpha_1 ( \begin{pmatrix} a_1 & a_2 \\ a_3 & a_4 \end{pmatrix}
)^{-1}$$
and so there exists a unique injective group homomorphism
$\hat \alpha : SL(2, \CC) \ltimes \CC^2 \to \A (A_1)$ extending $\hat
\alpha_1$ and $\hat \alpha_2$. Note that $\hat
\alpha_2 ( \begin{pmatrix} b_1 \\ b_2 \end{pmatrix}
 )=e^{\mathrm{ad}({b_1q+b_2p})}$.
\end{enumerate}

\begin{remark}
\label{precizie}
The derivative of $\hat \alpha_1 : SL(2, \CC) \to \A (A_1)$
is $\mathrm{ad} \circ f_I$ (cf. Example \ref{ex-I}) in the following sense: if $z\in\s$, then $\Phi_z :\CC\to \A (A_1)$ defined by $\Phi_z (t)=\hat \alpha_1 (e^{t\; z})$ is a group homomorphism which exponentiates $ {\mathrm{ad}}(f_I(z))$ in the sense of Definition \ref{definitie}.
\end{remark}

\begin{enumerate}
\item [{\bf E10}] For  $i_1, \ldots, i_n\in\NN^*$ we define
$R(i_1,\ldots, i_n)$ to be the (simply-connected) Lie group whose underlying manifold is $\CC^{n+1}$ and whose group law is
\beqa
\label{inmultire}
(a_1,\dots,a_n,v).(a'_1,\dots,a'_n,v')=(a_1+a'_1e^{-vi_1},\dots,a_n+a'_ne^{-vi_n},v+v').
\eeqa
\noi
One can check that the Lie algebra of $R(i_1,\ldots, i_n)$  is isomorphic to $\rr (i_1,\dots,i_n)$ and, using \eqref{inmultire}, one shows that
 $\Phi: R (i_1,\ldots, i_n) \to \A (A_1)$ defined by
\beqa
\Phi ((a_1,\dots,a_n,v)) (p)&=&e^{-v}p,\nonumber \\
 \Phi ((a_1,\dots,a_n,v)) (q)&=&e^v (q+\sum_{k=1}^n  \frac{a_k}{(i_k-1) !}
 p^{i_k-1}).\nonumber
\eeqa
\noi
 is a group homomorphism with discrete kernel isomorphic to $\ZZ$. Note that $\Phi ((a_1,\dots,a_n,v))=e^{\mathrm{ad}(\frac{a_1}{i_1!}p^{i_1}+\dots+\frac{a_n}{i_n!}p^{i_n})}e^{\mathrm{ad}(vpq)}$ which means that $\Phi$ exponentiates the formulae of example {\bf E8}.

\item [{\bf E11}] 
We define  $\tilde L_n$ to be the (simply-connected) Lie group whose underlying manifold is $\CC^{n+2}$ and whose group law is
\beqa
\label{grup-ln}
 (a_1,\dots, a_n,t,v) . (a'_1,\dots,a'_n,t',v')= (a''_1,\dots,a''_n,t'',v'')
\eeqa
\noi where 
\beqa
\label{coef-ln}
a''_k&=&a_ke^{(n-k)v'}+\sum_{j=1}^{k-1} \frac{t^{k-j}}{(k-j)!}a'_je^{-(k-j)v'}+a'_k,\nonumber\\
t''&=&t'+te^{-v'},\nonumber\\
v''&=&v+v'.\nonumber
\eeqa
\noi
One checks that the Lie algebra of $\tilde L_n$
is isomorphic to $\tilde {\cal L}_n$, that $R(n-1,n-2,\dots,0)$ is a subgroup of $\tilde L_n$ by the inclusion:
$$(a_1,a_2,\dots,a_n,v)\mapsto (a_1 e^{(n-1)v}, a_2 e^{(n-2)v},\dots,a_n,0, v)$$
and that $\tilde L_n/\left\{ (0,\dots,0,a_n,0,0):\  a_n\in\CC \right\}$ is a simply-connected Lie group whose Lie algebra is isomorphic to $\tilde {\cal L}_n/\CC$.
 In fact one can extend the map $\Phi$ of {\bf E10} to  $\tilde \Phi : \tilde L_n \to \A (A_1)$  by setting $\tilde \Phi ( (a_1,\dots,a_n,t,v))=\Phi ((a_1 e^{(-n+1)v}, a_2 e^{(-n+2)v},\dots,a_n,v)) e^{-t\mathrm{ad}(q)}$. Explicitly, this gives
\beqa
\label{muie-michel}
\tilde \Phi ( (a_1,\dots,a_n,t,v))(p)&=&e^{-v}p+t,\nonumber \\
\tilde \Phi ( (a_1,\dots,a_n,t,v))(q)&=&e^v \left(q + \sum_{k=1}^{n-1} \frac{a_k e^{(-n+k)v}}{(n-k-1)!}p^{n-k-1}\right).
\eeqa
Using \eqref{grup-ln} and \eqref{muie-michel} one checks that $\tilde \Phi$ is a group homomorphism whose kernel is the subgroup $\left\{ (0,\dots,0,a_n,0,2\pi i k):\ k\in\ZZ, a_n\in\CC \right\}$, isomorphic to the direct product $\CC\times\ZZ$. 
Finally, note that the subgroup $L_{n}$ of $\tilde L_n$ defined by $v=0$ is simply-connected and that its Lie algebra is isomorphic to ${\cal L}_{n}$. 
The restriction of $\tilde \Phi$ to $L_{n}$ 
factors to give an injection of ${ L}_n/\CC$ in $\A (A_1)$ and ${ L}_n/\CC$ is a simply-connected Lie group whose Lie algebra is isomorphic to ${\cal L}_{n-1}$.
\end{enumerate}

\section{A family $\cal N$ of $\s$ realisations in the Weyl algebra}
\renewcommand{\theequation}{5.\arabic{equation}}   
\setcounter{equation}{0}

In this section we study some explicit examples of realisations of
$\s$ in 
$A_1$ first given by Joseph in \cite{joseph}.
We show that distinct members of this family are inequivalent under the
action of the group $\A (A_1)\times \A (\s)$.
From now on we write:
$$e_+=\begin{pmatrix} 0 & 1 \\ 0 & 0 \end{pmatrix},
e_-=\begin{pmatrix} 0 & 0 \\ 1 & 0 \end{pmatrix}, e_0=\begin{pmatrix}
1 & 0 \\ 0 & -1 \end{pmatrix}$$
 and the images of $e_+, e_-,e_0$ under the natural inclusion 
$\s\subset {\cal U}(\s)$ will be  respectively denoted by $x,y,h$.

\begin{definition}
The group ${\mathrm{Aut}}(A_1)\times  {\mathrm{Aut}}({\cal
U}(\s))$ acts on $
A_1^{\s}$ by:
\beqa
(\alpha, w)\cdot f= \alpha\circ f\circ w^{-1} \vert_{\s}
\eeqa
where $\alpha \in \A (A_1)$, $w\in \A({\cal U}(\s))$ and $f \in  A_1^\s$.
Here ${\cal U}(\s )$ denotes the universal enveloping algebra of $\s$,
and we consider  $\s$ as embedded in ${\cal U} (\s)$.
\end{definition}

\begin{remark}
The group $\A(\s)$ is naturally included in the group $\A({\cal
U} (\s))$. This inclusion is
strict, see for example the equations (\ref{nou}). 
\end{remark}

The set $A_1^{\mathfrak
sl(2)}$ is in bijection with  the set of $\s$ triplets.
 If $X,Y,H$ is an $sl(2)$ triplet, $f:{\mathfrak {sl}(2)}\to A_1$ given by
$f(e_+)=X, f(e_-)=Y, f(e_0)=H$ is a
Lie algebra homomorphism and conversely, $f(e_+), f(e_-),
f(e_0)$ is an $\s$ triplet if
$f:\s\to A_1$ is a Lie algebra homomorphism.

\begin{example}
\label{ex-I} 
From {\bf E1} we know that
\begin{eqnarray}
X = -\frac 12  q^2, \
Y = \frac 12 p^2, \
H = \frac 12 (pq+qp),  \nonumber
\end{eqnarray}
form an $\mathfrak {sl}(2)$ triplet. 
We denote by  $f_I:\s\to A_1$  the
corresponding Lie algebra homomorphism. 
\end{example}

\noi
This realisation of $\s$  has the following properties:
\begin{enumerate}
\item [(i)] the set of eigenvalues of
${\mathrm{ad}}({f_{I}(e_0)})$ is $\ZZ$; 
\item [(ii)] $f_I(n)\in\Delta_1$  for any nilpotent $n\in\s\setminus \{ 0 \}$;
\item [(iii)] $f_I(s)\in\Delta_3$ for any semi-simple $s\in\s\setminus \{ 0 \}$.
\end{enumerate}
\noi
Property (i) follows from the equation ${\mathrm{ad}}({f_{I}(e_0)})(p^i q^j)=(j-i) p^i
q^j$; (ii) and (iii) follow from Proposition \ref{iar-dixmier}
and Remark \ref{simplificari} since $f_I(\s)=W_2$.

\begin{remark}
\label{operatori}
In the standard representation of $A_1$, the  operators $X,Y,H$ of Example
\ref{ex-I} are represented by the differential operators:
\beqa
X = - \frac 12  x^2, \
Y =  \frac 12  \frac {d^2}{dx^2}, \
H = x\frac{d}{dx}+\frac 12. \nonumber
\eeqa
\noi
It is interesting to notice that one can obtain $\s$ triplets
represented by differential operators of arbitrary order by applying
appropriate automorphisms of $A_1$ to this example.
\end{remark}

\begin{example}
\label{ex-IIA}
(\cite{joseph}) For  $b\in{\mathbb C}$,  the three elements of $A_1$
\beqa
X = (b + pq) q, \
Y = - p, \
H = 2pq + b, \nonumber
\eeqa
form  an $\mathfrak sl(2)$ triplet. We  denote by $f_{II}^{b}:\s\to A_1$ the
corresponding Lie algebra homomorphism. 
\end{example}

\noi
This realisation of $\s$ is fundamentally  different from Example \ref{ex-I} with respect to each of the three properties above:
\begin{enumerate}
\item [(i)] the set of eigenvalues of
${\mathrm{ad}}({f_{II}^{b}(e_0)})$ is $2\ZZ$; 
\item [(ii)]
there exists a nilpotent $n\in\s$ such that
$f_{II}^{b}(n)\in\Delta_2$;
\item [(iii)]
 there exists a semi-simple $s\in\s$ such that
$f_{II}^{b}(s)\in\Delta_4$. 
\end{enumerate}

\noi
Property (i) is obvious and properties (ii) and (iii) are consequences of the

\begin{lemma}
\label{alta}
$ \lambda X + \mu Y + \nu H \in \Delta_1 \cup \Delta_3$ iff $\lambda =0$.
\end{lemma}
{\bf Proof}: ($\Rightarrow$): One shows by  induction that 
$$\mathrm{ad}^n({\lambda X + \mu Y + \nu H}) (q) = n! \lambda^n a^n q^{n+1}+ h_n(q)$$
where $h_n(q)$ is
a polynomial in $q$ of degree at most $n$. It then follows that if $\lambda \ne 0$, $\big(
\mathrm{ad}^n({\lambda X + \mu Y + \nu H})(q)\big)_{n\in\NN}$ spans an infinite-dimensional vector space and
hence $\lambda X + \mu Y + \nu H\notin \Delta_1\cup \Delta_3$ by Corollary $6.6$ of \cite{Dixmier}. 

($\Leftarrow$): If $\lambda = 0$ the result follows from Proposition
\ref{iar-dixmier}. QED

\begin{remark}
\label{ex-IIB}
(\cite{joseph}) For  $b\in{\mathbb C}$,  the three elements of $A_1$
\beqa
X = - q, \
Y =  p (b +  pq),\
H = 2pq +b\nonumber
\eeqa
\noi
form an $\s$ triplet. The corresponding Lie homomorphism is easily seen to be 
$(\alpha,\tau)\cdot f_{II}^{-(b+2)} $
where $\tau\in\A(\s)$ and $\alpha\in \A(A_1)$ are given by $\tau(e_+)=e_-,\
\tau(e_-)=e_+$   and $\alpha (p)=q,\ \alpha(q)=-p$.
However, note that there does not exist $\beta\in\A (A_1)$ and $b'\in\CC$ such that $(\beta,Id)\cdot f_{II}^{b'}=(\alpha,\tau)\cdot f_{II}^{-(b+2)} $ 
since $(\alpha,\tau)\cdot  f_{II}^{-(b+2)}(e_-)\in\Delta_2$ and $f_{II}^{b'}(e_-)\in\Delta_1$.
\end{remark}

\begin{definition}
\label{familia}
Define   ${\cal N}\subseteq A_1^\g$ by
 $ {\cal N}=\{f_I,
f_{II}^{b}:  b\in\CC \}. $
\end{definition}

\subsection{Inequivalence of elements of $\cal N$ under $\A (A_1)\times \A(\s)$}

In \cite{joseph}  Joseph showed that $A_1^\s$ is a disjoint union $S_1\cup S_2$ 
where $f\in S_1$ (resp. $f\in S_2$) iff the set of eigenvalues of ${\rm ad} (f(e_0))$ is $2\ZZ$ (resp. $\ZZ$). 
Furthermore, he subdivided $S_1$ (resp. $S_2$) into a
disjoint union $S_{11}\cup S_{12}\cup \dots \cup S_{1\infty}$ (resp. $S_{21}\cup S_{22}\cup \dots $), 
showed that for  $1<r<\infty$ (resp. $ 1\le r<\infty$) the $S_{1r}$  (resp. $S_{2r}$) are stable under the 
action of $\A (A_1)\times \A ({\cal U} (\s))$ and that  $f_I\in S_{21}$ and $f_{II}^b\in S_{11}$. This means in particular that
$f_I$ and $f_{II}^b$ are inequivalent under $\A (A_1)\times \A  (\s)$. We will prove that $f_{II}^b$ and $f_{II}^{b'}$ are 
inequivalent under $\A (A_1)\times \A  (\s)$ if $b\ne b'$
so that distinct elements of $\cal N$ are inequivalent under the action of this group. It is not known  whether 
distinct elements of $\cal N$ are inequivalent under the action of $\A (A_1) \times \A({\cal U}(\s))$.

\begin{definition}
\label{casimir}
Let $f\in A_1^{{\mathfrak sl(2)}}$, let ${\hat f} : {\cal U} (\s)\to A_1$ be
the natural extension of $f$ to the universal enveloping algebra and let 
$Q_f={\hat f}(Q)$
be the image by $\hat f$ of the Casimir operator $Q=\frac 12 h^2 +xy+yx$ of $\s$.
\end{definition}

\noi
If  $w\in \A({\cal U}(\s))$ and $\alpha\in \A(A_1)$, we have 
$ Q_{{\hat f} \circ w \vert_\s}=({\hat f}\circ w) (Q)$ and $ Q_{\alpha\circ f}=Q_f$ (cf \cite{arnal} or \cite{dix}).
 Calculation shows that $Q_{f_{II}^{b}}=b(\frac 12 b +1)$ and thus:

\begin{proposition}
\label{cenume}
Let $b,b'\in  {\mathbb C}$ be such that $b'\ne b$ and $b'\ne -b-2$.
 Then 
 $ f_{II}^{b'}$ is not equivalent to $ f_{II}^b$
 under    $\A (A_1) \times \A({\cal U}(\s))$.
\end{proposition}

This  means that the only elements of $\cal N$ which can be equivalent
under the action of  $\A (A_1) \times \A({\cal U}(\s))$, {\it a fortiori} under the action $\A (A_1) \times\A(\s)$,  
are $f_{II}^{b}$ and  $f_{II}^{-b-2}$.  

\begin{proposition}
\label{memeCasimir}
If $b\ne -1$, $f_{II}^{b}$ is not equivalent to $f_{II}^{-b-2}$ under  $\A (A_1)\times \A (\s)$.
\end{proposition}
\textbf{Proof:} 
The $\s$ triplets corresponding to  $f_{II}^{b}$ and
$f_{II}^{-b-2}$ are
\beqa
\label{tablou222}
\begin{array}{ll}
X=(b+pq)q \quad &  X'=(-b-2+pq)q \cr
Y=-p & Y'=-p \cr
H=2pq+b & H'=2pq-b-2.
\end{array}
\eeqa
Suppose there exists $(\alpha,w)\in \A (A_1)\times \A (\s)$ such that $f_{II}^{b}=(\alpha,w)\cdot f_{II}^{-b-2}$,
i.e., such that $\alpha\circ \cdot f_{II}^{-b-2}=f_{II}^{b}\circ w$. Then writing $p'=\alpha(p), q'=\alpha(q)$ and
$w={\rm Ad}(g)$, this gives (see the Appendix) 
\beqa
\label{tablou22}
\begin{array}{ll}
(-b-2+p'q')q'&= a_1^2 (b+pq)q  - a_3^2
(-p) - a_1a_3 (2pq+b) \cr
\hskip 2.2truecm -p'&=-a_2^2 (b+pq)q  + a_4^2
(-p) + a_2a_4(2pq+b)  \cr
\hskip 0.5truecm 2p'q'-b-2&=-2 a_1a_2(b+pq)q  + 2a_3a_4 (-p) + (a_1 a_4 +a_2 a_3) (2pq+b). 
\end{array}
\eeqa
\noi
Substituting the second equation in the third equation, we obtain
\beqa
\label{identificare22}
-2[-a_2^2 (b+pq)q  + a_4^2
(-p) + a_2a_4(2pq+b)]q'-b-2 \nonumber \\
= -2 a_1a_2(b+pq)q  + 2a_3a_4 (-p) + (a_1 a_4 +a_2 a_3) (2pq+b). 
\eeqa
\noi
If $a_2\ne 0$, the expansion of $q'$ in the standard basis
can  have only a  constant term, otherwise the term $2a_2^2 (pq)q q'$ on the
LHS contains terms which are not present in the RHS; but then
$[p',q']=0$ which is a contradiction and hence $a_2 = 0$. 
Since $ a_1a_4-a_2a_3=1$ this implies that  $a_1a_4=1$.
Equation \ref{identificare22} now reduces to
\beqa
\label{identificare222}
-2[ a_4^2 (-p)]q'-b-2= 2a_3a_4 (-p) + (2pq+b)
\eeqa
\noi
and equating the constant term on both sides of (\ref{identificare222}) gives $b=-1$, a contradiction.
\begin{proposition}
\label{largire1}
 Distinct elements  of ${\cal N}$ are not equivalent
under the action of $\A (A_1)\times \A (\s)$.
\end{proposition}
{\bf Proof}: Immediate from  Corollary \ref{cenume}
and Proposition \ref{memeCasimir}. QED

\section{The orbit of $\cal N$ under $\A(A_1)\times\A(\s)$}
\renewcommand{\theequation}{6.\arabic{equation}}   
\setcounter{equation}{0}

In this section we  give various characterisations of the orbit of $\cal N$ under $\A(A_1)\times\A(\s)$ in terms of the Dixmier partition and also in terms of exponentiation. 
We then calculate the isotropy groups of each of the elements  $\cal N$ and
finally, for the sake of completeness, we give an explicit example of an element in $A_1^\s$ which is not in the orbit of ${\cal N}$.

\subsection{Characterisations in terms of the Dixmier partition}

\begin{definition}
\label{partitii}
\beqa
{\cal D}&=&\{ f\in A_1^\s : \exists (\alpha,
   w)\in\A(A_1)\times\A(\s) \mbox{ s.t } (\alpha, w)\cdot f\in {\cal N} \} \nonumber \\
{\cal D}_3&=&\{ f\in A_1^\s : \exists z\in \s\setminus\{0\} \mbox{ s.t }
f(z)\in\Delta_3\} \nonumber \\
{\cal D}^\prime_3&=&\{ f\in A_1^\s : \exists z\in \s\setminus\{0\} \mbox{ 
s.t }
\mathrm{ad}({f(z)}) \mbox{ has an eigenvector in } \Delta_3 \}
 \nonumber \\
{\cal D}_1&=&\{ f\in A_1^\s : \exists z\in \s\setminus\{0\} \mbox{ s.t }
f(z)\in\Delta_1\} \nonumber \\
{\cal D}^\prime_1&=&\{ f\in A_1^\s : \exists z\in \s\setminus\{0\}  \mbox{ 
s.t }
\mathrm{ad}({f(z)}) \mbox{ has an eigenvector in } \Delta_1 \}
 \nonumber \\
{\cal E}&=&\{ f\in A_1^\s : \exists z\in\s\setminus\{0\}  \mbox{ s.t} \ {\mathrm{ad}}f(z) \mbox{ can be exponentiated } \}
\nonumber
  \eeqa
\end{definition}

\noi
We now show that  the above  sets are the same.
This means in particular  that ${\cal N}$ is a set of
normal forms for the action of $\A (A_1)\times \A (\s)$ 
on ${\cal D}_3, {\cal D}^\prime_3, {\cal D}_1,{\cal D}_1'$ and ${\cal E}$.

\begin{theorem}
\label{TEOREMA}
${\cal D}={\cal D}_3={\cal D}^\prime_3={\cal D}_1={\cal D}_1'={\cal E}$.
\end{theorem}
\textbf{Proof:} First,  note that ${\cal D}_3,{\cal D}^\prime_3,{\cal
D}_1,{\cal D}^\prime_1$ and  ${\cal E}$  are stable  under the action of  $\A (A_1)\times \A (\s)$
and that the inclusions ${\cal D}_3\subseteq{\cal D}^\prime_3$ and ${\cal D}_1\subseteq{\cal D}^\prime_1$
are obvious. We have already seen that ${\cal N}\subseteq {\cal D}_3$, 
 ${\cal N}\subseteq {\cal D}_1$ and ${\cal N}\subseteq {\cal E}$ 
(cf  Examples \ref{ex-I}, \ref{ex-IIA} and Proposition \ref{caracterizare}) and  hence, to
prove the theorem it will be sufficient to show that  ${\cal D}_3\subseteq
{\cal D}$ ,  ${\cal D}'_3\subseteq {\cal D}_3$, ${\cal D}^\prime_1\subseteq {\cal D}_3$ and ${\cal E}\subseteq {\cal D}_1$.

\medskip

\noi
{\bf{$\mathbf{{\cal D}_3\subseteq{\cal D}}$}}: Let $f\in {\cal D}_3$.
By hypothesis there exists  $z\in\s$  such
that $f(z)\in \Delta_3$ and $z$ must be semi-simple by Remark \ref{simplificari}.
By rescaling we can always suppose that
the eigenvalues (in $\s$) of $\mathrm{ad}(z)$ are $-2,0$ and $2$
and then there exists $w\in\A(\s)$ such that $w^{-1}(e_0)=z$.
By Theorem $9.2$ of \cite{Dixmier}, there exist $\alpha\in\A
 (A_1), \mu \in \CC^*$ and $ \nu\in\CC$ such that 
$(\alpha, w)\cdot f(e_0)=\alpha \circ f\circ w^{-1} (e_0)=\mu pq + \nu$. 
By Lemma $2.4$ of \cite{joseph} and Remark \ref{ex-IIB}, there exists
$(\alpha', w')\in\A (A_1)\times \A (\s)$ such
$(\alpha', w')\cdot((\alpha, w)\cdot f)\in {\cal N}$ and hence $f\in{\cal D}$.

\medskip
\noi
{\bf {${\cal D}^\prime_3\subseteq{\cal D}_3$}}: We need the following lemma

\begin{itemize} \item []
\begin{lemma}
Let $S\in \Delta_3$. Then 
\beqa
\label{delta5}
\begin{array}{llll}
(i)&C(S)&\subset& \CC \cup \Delta_3 \cup \Delta_5;\\
(ii)& C(S)\cap \Delta_3 &=&\{ \mu S + \nu : \mu\in\CC^*, \nu\in\CC\}.
\end{array}
\eeqa
\end{lemma}
\textbf{Proof:} Without loss of generality we can suppose that $S=pq$ since
any element of $\Delta_3$ is equivalent under $\A (A_1)$ to $\mu' pq + \nu'$
(Theorem $9.2$ of \cite{Dixmier}) and since $C(\mu' pq + \nu')= C(pq)$.
First note that $C(pq)=\CC [pq]$ (see Proposition $5.3$ of
\cite{Dixmier}). Let $Z=a_k (pq)^k + \dots+a_0$ be a polynomial of degree $k$ in
$pq$. Then a simple induction shows that 
$$[Z,p^m q^n]= k (n-m) a_kp^{m+k-1} q^{n+k-1} + \mbox{ terms of lower degree in
$p$}. $$
From this it follows that if $k>1$ the only eigenvalue of ${\rm ad}(Z)$ is $0$
and so $D(Z)= C(Z)$. By iteration of this formula it also follows that that if
$k>1$, ${\rm Ker} \, {\rm ad}^\ell (Z) = {\rm Ker} \, {\rm ad}(Z) $ and so $N(Z)=C(Z)$. 
This means that if $k>1$, $Z\in\Delta_5$ (see Theorem $2.2$).  It is clear
that  $Z\in \Delta_3$ if $k=1$ and the lemma is proved. QED

\end{itemize}

Let  $f \in {\cal D}^\prime_3$.  There exists $z\in \s \setminus\{0\}$,
$S\in \Delta_3$ and $\lambda\in \CC$ such that $[f(z),S]=\lambda S$. 
Since $S\in\Delta_3$, $f(z)$
is a sum of eigenvectors of ${\rm ad}(S)$ from which it follows  that $f(z)\in
C(S)$. By Lemma \ref{delta5}, $f(z)\in \Delta_3 \cup \Delta_5$; but $f(z)\notin
\Delta_5$ (see Remark $3.2$) and so $f(z)\in \Delta_3$ and  $f\in{\cal D}_3$.

\medskip
\noi
{\bf{${\cal D}'_1\subseteq{\cal D}_3$}}:  We need the following two lemmas

\begin{itemize} \item []
\begin{lemma}
\label{lema2}
Let $a\in A_1$, $\mu \in \mathbb C^*$ be such that
$[a,p]=-\mu p $. There exists  $\alpha_3\in \A (A_1)$ and $a_0\in\CC$
such that $\alpha_3^{-1}(a)=\mu pq+a_0$.
\end{lemma}
\textbf{Proof:}
Let $a=\sum_{ij}h_{ij} p^iq^j$. 
Since 
\beq
\label{condHmax}
[a,p]=-\mu p \ \Leftrightarrow\ [a-\mu pq,p]=0,
\eeq
\noi and $[p,p^iq^j]=jp^iq^{j-1}$, one has
\beq
\label{muie}
a=\mu pq + \sum_{i=0}^N a_i p^i,  
\eeq
where $a_i \in \mathbb C$ and $N\in \mathbb N$.
One can then write 
$$a=\mu p(q+\sum_{i=1}^N \frac{a_i}{\mu}p^{i-1})+a_0. $$

But $[p,q']=1$ where $q'=q+\sum_{i=1}^N \frac{a_i}{\mu} p^{i-1}$ and it is easy to see that the
 homomorphism $\alpha_3:A_1\rightarrow A_1$
given by $\alpha_3(p)=p,\ \alpha_3(q)=q'$ is invertible.
Hence
$$\alpha_3^{-1}(a)=\mu pq + a_0. $$
{\flushright QED}

\begin{lemma}
\label{lema1}
Let $a\in A_1$, $\lambda\in\CC^*$ and let $g(p)= \sum_{k=0}^n
b_k p^k$ be a polynomial of degree $n$ in $p$.
If 
$[a, g(p)]=\lambda g(p)$ then there
exists $\alpha_2\in \A(A_1)$ such that 
$[\alpha_2^{-1}(a),
p]=\frac{\lambda}{n} p$.
\end{lemma}
\noi
\textbf{Proof:} Since $<p^iq^j>$ is a basis of $A_1$, we can write
$$[a,p^k]=\sum_{i=0}^{N_k} f_{i,k}(p)\ q^i, $$
where $k\in {\mathbb N}^*, N_k\in \mathbb N$ and $f_{N_k,k}$ is a non-zero polynomial of degree $M_k$ in
$p$. Let $\alpha_k\ne 0$ be the coefficient of $p^{M_k}q^{N_k}$ in $[a,p^k]$.
A straightforward induction argument shows that $N_k=N_1,\ M_k=M_1+k-1$ and $\alpha_k=k\alpha_1$.

Thus, one has
$$[a,\sum_{k=0}^n b_k p^k]=\sum_{k=0}^n \left(\sum_{i=0}^{N_1} b_k f_{i,k}(p) q^i \right). $$
 But since $[a, g(p)]=\lambda g(p)$ by hypothesis, this means that
 $b_n \alpha_n=b_n n \alpha_1= b_n \lambda$, 
 $N_1=0,\ M_n=n$ and so $\alpha_1=\frac{\lambda}{n}$ and $M_1=1$. Therefore
$$[a,p]=\frac{\lambda}{n} p + \nu$$
for some constant $\nu \in \mathbb C$.
But $[\frac{\lambda}{n} p + \nu, \frac{n}{\lambda} q]=1$ 
and so there exists an unique automorphism $\alpha_2$ of $A_1$
such that
$\alpha_2(p)=\frac{\lambda}{n} p + \nu$ and $ \alpha_2(q)=\frac{n}{\lambda} q$.
Hence
$$[\alpha_2^{-1}(a),p)]=\frac{\lambda}{n} p. $$
{\flushright {QED}}
\end{itemize}

\noindent We now prove that ${\cal D}'_1\subseteq {\cal
D}_3$. Let $f\in {\cal D}'_1$.
By hypothesis there exist  $z\in\s$,
$N\in \Delta_1$ and $\lambda\in\CC$ such that
$$ [f(z), N]=\lambda N.$$
By Theorem $9.1$ of 
\cite{Dixmier} there exist $\alpha_1\in\A(A_1)$ such that
$\alpha_1(N)\in\CC[p]$. If  $\lambda\ne 0$,
$\alpha_1\circ f(z) $ satisfies the hypothesis of Lemma \ref{lema1}. 
and, by  Lemmas \ref{lema1} and  \ref{lema2},
 there exist
$\alpha_2, \alpha_3\in\A (A_1)$, $\mu\in \CC^*$ and $a_0\in\CC$ such that $\alpha_3^{-1}\circ
\alpha_2^{-1}\circ\alpha_1\circ f(z) = \mu pq + a_0$. 
Since $\mu pq + a_0$ is in $\Delta_3$,  this means that $\alpha_3^{-1}\circ
\alpha_2^{-1}\circ\alpha_1\circ f\in {\cal D}_3$ and hence $f\in {\cal D}_3$.
To complete the proof we show that the $\lambda=0$ case reduces to the 
$\lambda\ne 0$ case as follows.
If $\lambda=0$, then $f(z)\in C(N)$. But $C(N)\subseteq\Delta_1$ (by 
Theorem $9.1$ of 
\cite{Dixmier}) so that $f(z)\in\Delta_1$ and, by Remark \ref{simplificari},
 $z$ is nilpotent. There exists $s\in\s$ semi-simple such that
 $[s,z]=2z$ and then $[f(s),f(z)]=2f(z)$. 

\medskip
\noi
{\bf ${\cal E}\subseteq {\cal D}_1$}:
If $f\in{\cal E}$ then by  the above and Proposition \ref{caracterizare},  $f\in{\cal D}_3 \cup {\cal D}_1={\cal D}_1$. QED

\begin{corollary}
Let $f\in A_1^\s$. The following are equivalent:
\begin{enumerate}
\item
$f(\s)\subseteq \Delta_1 \cup \Delta_3$.
\item
For all $z\in \s$, $f(z)$ can be exponentiated.
\item 
There exists $(\alpha, w)\in \A (A_1)\times \A(\s)$ such that $(\alpha, w).f=f_I$.
\end{enumerate}
\end{corollary}

\subsection{Isotropy groups}

Recall that if a group $G$ acts on a set $X$,  the
isotropy of $x\in X$ in $G$ is by definition
the subgroup $\{ g\in G : \ g.x=x \}.$
In this subsection we calculate the isotropy of $f_I$ and $f_{II}^{b}$ under the action of the group $\A(A_1)\times\A(\s)$.

In order to  calculate the isotropy of $f_I$  we need the following lemma. As usual, 
 $\mathrm{Ad}:SL(2,\CC)\to
\A(\s)$ is defined by 
 $\mathrm{Ad} (g)
 (x)=gxg^{-1}$, where $g\in SL(2,\CC)$ and  $x\in\s$.

\begin{lemma}
\label{lema-iso-I}
There exists a group homomorphism $\hat \alpha_1: SL(2,\CC) \to \A(A_1)$
such that $f_{I}\circ \mathrm{Ad} (g) = \hat \alpha_1 (g) \circ
f_{I}$ for all $g\in SL(2,\CC)$. Furthermore, 
 $\hat \alpha_1$ does not factor through $\mathrm{Ad}:SL(2,\CC)\to
\A(\s)$.
\end{lemma}
{\bf Proof:} 
Define  $\hat \alpha_1 : SL(2,\CC) \to \A(A_1)$ by (see {\bf E9}, subsection $4.4$) 
$$\hat \alpha_1
(\begin{pmatrix} a_1 & a_2 \cr a_3 & a_4 \end{pmatrix}) (p)=a_2q+ a_4
p,\ \hat \alpha_1 (\begin{pmatrix} a_1 & a_2 \cr a_3 & a_4  \end{pmatrix}) (q)=
a_1 q + a_3 p.$$
The formula   $f_{I}\circ \mathrm{Ad} (g) = \hat \alpha_1 (g) \circ
f_{I}$ follows by a straightforward calculation and
$\hat \alpha_1$ cannot factor through $\mathrm{Ad}:SL(2,\CC)\to
\A(\s)$ since $\hat \alpha_1 (-Id)\ne \hat \alpha_1 (Id)$. QED

\begin{proposition}
\label{iso-I}
Let  ${\cal I}_{f_I}$ be the isotropy of  $f_I$ in 
$\A(A_1)\times\A(\s)$. Then
$$
{\cal I}_{f_I}=
\big\{ (\hat \alpha_1(g), \mathrm{Ad} (g))\ : \ g\in SL(2,\CC)\big\}.
$$ 
In particular, ${\cal I}_{f_I}$ is isomorphic to $SL(2,\CC)$.
\end{proposition}
{\bf Proof:} The inclusion  $\big\{ (\hat \alpha_1(g), \mathrm{Ad} (g))\ : \ g\in SL(2,\CC)\big\}\subseteq {\cal I}_{f_I}$ 
follows immediately from the lemma. To prove inclusion in the opposite sense, let $(\alpha, w)\in {\cal I}_{f_I}$ 
and choose $g\in SL(2,\CC)$ such that $ w=\mathrm{Ad} (g)$. 
By definition 
$$f_I\circ \omega=\alpha\circ f_I$$
 and, since  by the previous
lemma 
$f_I\circ
 \omega=\hat \alpha_1 (g)\circ f_I$, we get  $\alpha^{-1}\circ\hat
 \alpha_1(g)\circ f_I= f_I$.
Using the explicit formulae for $f_I$ (cf Example \ref{ex-I}) this implies
 that $q'^2=q^2$, $p'^2=p^2$ where we have written
 $p'=\alpha^{-1}\circ \hat \alpha_1(g) (p)$ and
 $q'=\alpha^{-1}\circ \hat \alpha_1(g)(q)$.
Hence $q'\in C(q^2)=\CC[q]$
 and $p'\in C(p^2)=\CC[p]$. It is then easy to see that either $q'=q$ and
  $p'=p$, or $q'=-q$ and  $p'=-p$. Thus either $\alpha^{-1}\circ\hat \alpha_1(g)=1$ or
 $\alpha^{-1}\circ\hat \alpha_1(g)=\hat \alpha_1(-Id)$  which means either
$\alpha=\hat \alpha_1(g)$ or $\alpha= \hat \alpha_1(-g)$. Since $w=\mathrm{Ad} (g)=\mathrm{Ad} (-g)$, 
this completes the proof of the proposition. QED

\medskip

To calculate the isotropy of $f_{II}^{b}$ we need the following
definition and lemma:

\begin{definition}
\label{B}
Let $\hat B=\left\{\begin{pmatrix} a_1 & 0 \cr a_3 & \frac{1}{a_1} \end{pmatrix}\in SL(2,
\CC): a_1\in\CC^*,\ a_3\in\CC\right\}$ and let $B$ be the subgroup
$\mathrm{Ad}(\hat B)$ of $\A(\s)$.
\end{definition}

\begin{lemma}
\label{lema-iso-IIA}
There exists a  group homomorphism $\hat \beta: \hat B \to \A(A_1)$
such that for all $g\in \hat B$,
$f_{II}^{b}\circ \mathrm{Ad}(g )= \hat \beta (g) \circ
f_{II}^{b}$. Furthermore, $\hat \beta (g) = \hat \beta (-g)$
and
$\hat \beta$ factors through $\mathrm {Ad}: \hat B\to B$ to a group homomorphism $\beta:B\to
\A (A_1)$.

\end{lemma}
\begin{remark}
\label{altaremarca2}
The derivative of $\hat \beta: \hat B \to \A(A_1)$ is the restriction
of the Lie
homomorphism $\mathrm{ad}\circ f_{II}^{b} : \s \to \mathrm{Der} (A_1)$ to
the Lie algebra $\mathfrak b\subseteq\s$ of $\hat B$.
This makes sense since every element of $\mathrm{ad}\circ f_{II}^{b} ({\mathfrak b})$ is exponentiable. 
Note that the elements of $\mathrm{ad}\circ f_{II}^{b} \big(\s\setminus{\mathfrak b}\big)$ are not exponentiable by Lemmas \ref{alta} and  \ref{caracterizare}.
\end{remark}
{\bf Proof:} Define  $\hat \beta : \hat B \to \A(A_1)$ by 
$$\hat \beta (
\begin{pmatrix} a_1 & 0 \cr a_3 & \frac{1}{a_1} \end{pmatrix}) (p)=
\frac {1}{a_1^2} p,\ \hat \beta (
\begin{pmatrix} a_1 & 0 \cr a_3 & \frac{1}{a_1} \end{pmatrix}) (q)=
a_1^2 (q-\frac{a_3}{a_1}).$$
The result then follows by a straightforward calculation. QED

\begin{proposition}
\label{iso-IIA} 
Let ${\cal I}_{f_{II}^b}$ be the isotropy of  $f_{II}^{b}$ in 
$\A(A_1)\times\A(\s)$. Then 
$${\cal I}_{f_{II}^b}=\{ (\beta(w), w)\ : \ w\in B\}. $$ 
In particular, ${\cal I}_{f_{II}^b}$ is isomorphic to a Borel subgroup of $SO(3,\CC)$.
\end{proposition}
{\bf Proof:} The inclusion  $\{ (\beta(w), w)\ : \ w\in B\}\subseteq {\cal I}_{f_{II}^b}$ 
follows immediately from the lemma. To prove inclusion in the opposite sense, let $(\alpha, w)\in {\cal I}_{f_{II}^b}$ 
and choose $g\in SL(2,\CC)$ such that $ w=\mathrm{Ad} (g)$.  By definition,  we have
\beqa
\label{muie4}
 \alpha\circ f_{II}^{b}= f_{II}^{b}\circ w,
\eeqa
which is equivalent to
$$ 
\alpha \circ f_{II}^{b}(x)= f_{II}^{b} (gxg^{-1})\qquad\forall x\in\s.
$$
If $x\in  {\mathfrak b}=<e_0,e_{-}>$ 
then  $f_{II}^{b}(x)\in\Delta_1\cup \Delta_3$ (by Lemma
\ref{alta}),
$ \alpha \circ f_{II}^{b}(x)\in\Delta_1\cup \Delta_3$ (since the Dixmier partition is invariant under $\A(A_1))$ and hence
$f_{II}^{b} (gxg^{-1})\in \Delta_1 \cup \Delta_3$.
By Lemma \ref{alta}, this means that $gxg^{-1}\in {\mathfrak
b}$ and we have shown that $g  {\mathfrak
b} g^{-1}= {\mathfrak b}$.
Since $ {\mathfrak
b}$ is a Borel subalgebra of $\s$ and is the Lie algebra of $\hat B$, this implies that $g\in \hat B$. 

By the previous lemma and equation \eqref{muie4}, we have 
$$ 
\alpha  \circ f_{II}^{b}= f_{II}^{b}\circ
w=f_{II}^{b}\circ \mathrm{Ad}( g)= \hat \beta (g)
\circ f_{II}^{b} = \beta (w)\circ f_{II}^{b} 
$$
and hence 
$$
\beta (w)^{-1}\circ \alpha \circ f_{II}^{b} = f_{II}^{b}.
$$
Writing $p'=\beta (w)^{-1}\circ \alpha (p)$, $q'=\beta (w)^{-1}
\circ \alpha (q)$ and using the explicit formulae for $f_{II}^{b}$ (cf Example \ref{ex-IIA}),  this implies that $p'=p$ and
$2p'q'+b=2pq+b$. 
Hence  $0=2(pq'-pq)=2p(q'-q)$, which implies that $q'=q$ and $\alpha=\beta(w)$. QED

\subsection{Other examples}

One can construct elements of $A_1^\s$ which are not in the orbit
of  $\cal N$ under $\A(A_1)\times\A(\s)$
by letting
$\A ({\cal U}(\s))\setminus\A(\s)$  act on
 $\cal N$. In this section, for the sake of completeness, we
give an explicit example (see also page $127$ of \cite{joseph}).

Define $g\in A_1^\s$ by $g=f_{II}^{1}\circ w$, where 
$w=\exp(\mathrm{ad}({x^2}))\in\A ({\cal U}(\s))$ is given by:
\beqa
\label{nou}
w(x)&=&x,\nonumber\\
w(y)&=&y+hx+xh-4x^3,\nonumber\\
w(h)&=&h-4x^2.
\eeqa 
\noi
Explicitly, this gives
\beqa
\label{exotic}
g(e_+)&=&(1+pq)q\nonumber \\
g(e_-)&=&-p+4p^2q^3-4p^3q^6+12p^2q^5\nonumber \\
g(e_0)&=&2pq-4p^2q^4+1.
\eeqa

\begin{proposition}
\label{compliteness}
$g\notin\cal D$.
\end{proposition}
\textbf{Proof:} It is enough to prove that there does not
exist $(\alpha,w) \in \A(A_1)\times \A(\s)$ and $f\in \cal N$ such that $(\alpha,w).f=g $. 
First note that if $(\alpha,w).f=g$ then either
$f=f_{II}^{1}$ or $f=f_{II}^{-3}$ by 
 Proposition \ref{casimir} and the formula
$Q_{f_{II}^{b}}=\frac 12 b (b+1)$.
Suppose there exists  $(\alpha,w) \in
\A(A_1)\times \A(\s)$ such that $(\alpha,w).f_{II}^{b}=g$  with $b\in\{1,-3\}$.

The $\s$ triplets corresponding to $g$ and $f_{II}^{b}$
are
\beqa
\label{tablou-nou}
\begin{array}{ll}
X= (1+pq)q &  X'=(b+pq)q \cr
Y=-p+4p^2q^3-4p^3q^6+12p^2q^5  & Y'=-p \cr
H=2pq+1-4p^4q^2 & H'=2pq+b.
\end{array}
\eeqa
Writing
 $p'=\alpha(p),\ q'=\alpha(q)$ and  $w=
{\rm Ad}(\begin{pmatrix} a_1 & a_2 \cr a_3 & a_4 \end{pmatrix})$, this gives (see appendix)
\beqa
\label{tablou-nou2}
(b+p'q')q'&=&a_1^2X -a_3^2 Y - a_1a_3 H\nonumber \\
&=& a_1^2 ((1+pq)q)-a_3^2 (-p+4p^2q^3-4p^3q^6+12p^2q^5) -a_1a_3 (2pq+1-4p^4q^2) \nonumber \\
-p'&=& -a_2^2 X  +a_4^2 Y + a_2a_4 H \nonumber \\
&=& -a_2^2(1+pq)q  +a_4^2 (-p+4p^2q^3-4p^3q^6+12p^2q^5) + a_2a_4 (2pq+1-4p^4q^2)\nonumber  \\
2p'q'+b&=& -2a_1a_2 X + 2a_3a_4 Y + (a_1a_4+a_2a_3) H\nonumber \\
&=&-2a_1a_2 (1+pq)q + 2a_3a_4 (-p+4p^2q^3-4p^3q^6+12p^2q^5)  \nonumber\\
&&+ (a_1a_4+a_2a_3)( 2pq+1-4p^4q^2) .\nonumber
\eeqa
\noi
Substituting the second equation in the third equation, we obtain
\beqa
\label{master}
-2a_1a_2 (1+pq)q + 2a_3a_4 (-p+4p^2q^3-4p^3q^6+12p^2q^5) + (a_1a_4+a_2a_3)( 2pq+1-4p^4q^2)\nonumber\\
=-2\Bigl( -a_2^2(1+pq)q  +a_4^2 (-p+4p^2q^3-4p^3q^6+12p^2q^5) + a_2a_4 (2pq+1-4p^4q^2)\Bigr)q' +b.\nonumber
\eeqa

If $a_4\ne 0$ then the expansion of $q'$ in the standard basis
would consist of only a scalar term otherwise  $-2 a_4^2 p^3q^6 q'$ on the RHS
contains terms  which are not present on the LHS; 
but then $[p',q']=0$ which is a contradiction and hence $a_4=0$.

The equation above now reduces to
\beqa
\label{master2}
-2a_1a_2 (q+pq^2)  - (2pq+1-4p^4q^2)
= 2a_2^2(q+pq^2) q' +b.
\eeqa
\noindent Let $k$ be the highest power of $q$ appearing in the expansion of $q'$
in the standard basis. Then the highest power appearing in the
expansion of the RHS in the standard basis is $k+2$. Comparing with
the LHS gives $k=0$ and $q'=f(p)$ where $f$ is a polynomial in $p$ of degree  at most $3$.
This polynomial must in fact be of degree $3$ to provide the term $4p^4q^2$ on the LHS
but then this introduces a term in $p^3q$ on the RHS which is not present on
the LHS. In conclusion \eqref{master2} has no solutions and thus $g\notin
{\cal D}$.
QED

\begin{corollary}
$g(\s)\subseteq \Delta_2\cup\Delta_4$.
\end{corollary}
{\bf Proof:} Immediate from Theorem \ref{TEOREMA}. QED

\begin{remark}
If we set
$${\cal D}_1^{''}=\{ f\in A_1^\s : \exists n\in \s \mbox{ nilpotent
s.t }
{f(n)} \mbox{ is an eigenvector 
of some element in }\Delta_3\}$$
then $g\in {\cal D}''_1$
since $g(e_+)$ is an eigenvector of $pq$. Hence $ {\cal
D}''_1\ne  {\cal D}$.
\end{remark}

With respect to the Joseph decompositions 
$S_1=S_{11}\cup S_{12}\cup \dots \cup S_{1\infty}$ and $S_2=S_{21}\cup S_{22}\cup \dots $, 
one can show that  $g\in S_{1\infty}$ (see \cite{joseph} page $127$ and the Appendix).
It is not known whether  $S_{ij}$ is non-empty when $j>1$ and so 
it seems reasonable to conjecture that in fact, all elements of $A_1^\s$ are 
obtained from $\cal N$ by the action of the group $\A (A_1)\times \A ({\cal U} (\s))$.

\appendix

\section{Appendix}
\renewcommand{\theequation}{A.\arabic{equation}}   
\setcounter{equation}{0}

\begin{proposition}
\label{apendice2}
Let  $g=\begin{pmatrix} a_1 & a_2 \cr a_3 & a_4 \end{pmatrix}\in SL(2, \CC)$. Define $w\in \A (\s)$ by $w(z)=gzg^{-1}$ for any $z\in\s$. 
Then
\beqa
\label{auto-sl2-2}
w(e_+)&=&\begin{pmatrix} -a_1a_3 & a_1^2 \cr -a_3^2 & a_1a_3 \end{pmatrix} \nonumber \\
w(e_-)&=&\begin{pmatrix} a_2a_4 & -a_2^2 \cr a_4^2 & -a_2a_4 \end{pmatrix} \nonumber\\
w(e_0)&=&\begin{pmatrix} a_1a_4+a_2a_3 & -2a_1a_2 \cr 2a_3a_4 & -(a_1a_4+a_2a_3) \end{pmatrix}, \nonumber
\eeqa
where $e_+=\begin{pmatrix} 0 & 1 \cr 0 & 0
\end{pmatrix}, e_-=\begin{pmatrix} 0 & 0 \cr 1 & 0 \end{pmatrix},
e_0=\begin{pmatrix} 1 & 0 \cr 0 & -1 \end{pmatrix}$.\\
\end{proposition}

\begin{proposition}
\label{S1infty}
Let $g\in A_1^\s$ be defined by equations \eqref{exotic}. Then $g\in S_{1\infty}$.
\end{proposition}
{\bf Proof:} By \cite{joseph}, $S_1=S_{11}\cup S_{12}\cup \dots \cup S_{1\infty}$ and $S_{1r}$  ($1<r<\infty$) are stable under 
$\A (A_1)\times \A ({\cal U} (\s))$. Hence to prove that $g\in S_{1\infty}$
it is enough to show that $g\notin S_{11}$. Set $H=g(e_0),\  X=g(e_+)$,  $Y=g(e_-)$ and,
 for $m\in 2\ZZ$, define $D\big(H,m\big)$ by
$$D\big(H,m\big)=\left\{ z\in A_1 :\ [H,z]=mz\right\}.$$
Then, by definition (see \cite{joseph}), $g\notin S_{11}$ iff
$D\big(H,2\big)=X \CC[H] \mbox{ and } D\big(H,-2\big)=Y \CC[H]$.
By Lemma $3.1$ of \cite{joseph}, there exists $y_2\in D\big(H,2\big)$ such that $D\big(H,2\big)=y_2 \CC[H]$. Since $X\in D\big(H,2\big)$ there exists a polynomial $a_n H^n+\dots+a_0$ such that
$$ X=y_2(a_n H^n+\dots+a_0),$$
which gives
$$ q+pq^2= y_2 \left(a_n (2pq-4p^2q^4+1)^n+\dots+a_0\right).$$
Comparing the highest power of $p$, one has $n=0$, $y_2=\frac{1}{a_0} X$ and thus $D\big(H,2\big)=X \CC[H]$.
A similar but slightly more complicated argument shows that $ D\big(H,-2\big)=Y \CC[H]$. QED

\end{document}